\documentclass[a4paper,fleqn]{cas-sc}

\usepackage[numbers]{natbib}
\usepackage{mathtools}
\usepackage{tikz}
\definecolor{fc}{HTML}{1E90FF}
\definecolor{h}{HTML}{228B22}
\definecolor{bias}{HTML}{87CEFA}
\definecolor{noise}{HTML}{8B008B}
\definecolor{conv}{HTML}{FFA500}
\definecolor{pool}{HTML}{B22222}
\definecolor{up}{HTML}{B22222}
\definecolor{view}{HTML}{FFFFFF}
\definecolor{bn}{HTML}{FFD700}
\tikzset{fc/.style={black,draw=black,fill=fc,rectangle,minimum height=0.5cm}}
\tikzset{h/.style={black,draw=black,fill=h,rectangle,minimum height=0.5cm}}
\tikzset{bias/.style={black,draw=black,fill=bias,rectangle,minimum height=0.5cm}}
\tikzset{noise/.style={black,draw=black,fill=noise,rectangle,minimum height=0.5cm}}
\tikzset{conv/.style={black,draw=black,fill=conv,rectangle,minimum height=0.5cm}}
\tikzset{pool/.style={black,draw=black,fill=pool,rectangle,minimum height=0.5cm}}
\tikzset{up/.style={black,draw=black,fill=up,rectangle,minimum height=0.5cm}}
\tikzset{view/.style={black,draw=black,fill=view,rectangle,minimum height=0.5cm}}
\tikzset{bn/.style={black,draw=black,fill=bn,rectangle,minimum height=0.5cm}}
\usetikzlibrary{decorations.pathreplacing,calligraphy}
\usetikzlibrary {arrows.meta}
\usepackage{subfig}

\def\tsc#1{\csdef{#1}{\textsc{\lowercase{#1}}\xspace}}
\tsc{WGM}
\tsc{QE}
\tsc{EP}
\tsc{PMS}
\tsc{BEC}
\tsc{DE}
\begin{document}
\let\WriteBookmarks\relax
\def\floatpagepagefraction{1}
\def\textpagefraction{.001}

\shorttitle{NSPOD}

\shortauthors{Francesc Levrero-Florencio et~al.}

\title [mode = title]{NSPOD: Accelerating Krylov solvers via DeepONet-learned POD subspaces\tnotemark[1]}                      
% \tnotemark[1,2]

\tnotetext[1]{G.E.K. is supported by the ONR Vannevar Bush Faculty Fellowship. We also acknowledge support from the DOE-MMICS SEA-CROGS DE-SC0023191 award, the Ansys Inc, and the DARPA-DIAL grant HR00112490484.}

% \tnotetext[2]{The second title footnote which is a longer text matter
%    to fill through the whole text width and overflow into
%    another line in the footnotes area of the first page.}

\author[1]{Francesc Levrero-Florencio}[]

\cormark[1]
\ead{flevrero@synopsys.com}
\credit{Conceptualization, Methodology, Software, Writing}

\affiliation[1]{organization={Synopsys UK Ltd},
    country={UK}}

\author[2]{Youngkyu Lee}[]
   
\cormark[1]
\ead{youngkyu_lee@brown.edu}
\credit{Conceptualization, Methodology, Software, Writing}

\affiliation[2]{organization={Division of Applied Mathematics, Brown University},
    state={Rhode Island},
    country={USA}}

\author[3]{Jay Pathak}[]
\cormark[2]
\ead{jpathak@synopsys.com}
\credit{Conceptualization, Writing}

\affiliation[3]{organization={Synopsys Inc},
    state={California},
    country={USA}}

\author[2]{George Em Karniadakis}[]
\cormark[2]
\ead{george_karniadakis@brown.edu}
\credit{Conceptualization, Methodology, Writing}  

\cortext[cor1]{These two authors contributed equally.}
\cortext[cor2]{Corresponding authors}

% \fntext[fn1]{This is the first author footnote. but is common to third
% \fntext[fn2]{Another author footnote, this is a very long footnote and
%   it should be a really long footnote. But this footnote is not yet
%   sufficiently long enough to make two lines of footnote text.}

% \nonumnote{This note has no numbers. In this work we demonstrate $a_b$
%   the formation Y\_1 of a new type of polariton on the interface
%   between a cuprous oxide slab and a polystyrene micro-sphere placed

\begin{abstract}
The convergence of Krylov-based linear iterative solvers applied to parametric partial differential equations (PDEs) is often highly sensitive to the domain, its discretization, the location/values of the applied Dirichlet/Neumann boundary conditions, body forces and material properties, among others. We have previously introduced hybridization of classical linear iterative solvers with neural operators for specific geometries, but they tend to not perform well on geometries not previously seen during training \cite{Zhang2024}. We partially addressed this challenge by introducing the deep operator network Geo-DeepONet and hybridizing it with Krylov-based iterative linear solvers \cite{Lee20251}, which, despite learning effectively across arbitrary unstructured meshes without requiring retraining, led to only modest reductions in iterations compared to \textit{state-of-the-art} preconditioners. In this study we introduce Neural Subspace Proper Orthogonal Decomposition (NSPOD), a multigrid-like deep operator network-based preconditioner which can dramatically reduce the number of iterations needed for convergence in Krylov-based linear iterative solvers, even when compared to \textit{state-of-the-art} methods such as algebraic multigrid preconditioners. We demonstrate its efficiency via numerical experiments on a linearized version of solid mechanics PDEs applied to unstructured domains obtained from complex CAD geometries. We expect that the findings in this study lead to more efficient hybrid preconditioners that can match, or possibly even surpass, the convergence properties of the current \textit{gold standard} preconditioning methods for solid mechanics PDEs.
\end{abstract}

% \includegraphics{figs/grabs.pdf}

%\item Research highlights item 1
%\item Research highlights item 2
%\item Research highlights item 3

% Keywords

%%%%%%%%%%%%%%%%%%%%%%%%%%%%%%%%%%%%%%%%%%%

\maketitle

%%%%%%%%%%%%%%%%%%%%%%%%%%%%%%%%%%%%%%%%%%%

\section{Introduction}
\noindent Numerical solutions of parametric elliptic partial differential equations (PDEs) are essential in real-world engineering applications. These solutions are obtained by applying a suitable numerical method to the corresponding PDEs, with the most commonly used numerical method being the finite element method (FEM), and then solving a potentially large system of linear equations \cite{Bathe, Ciarlet}. Depending on how large and complex the considered domain is, the number of nodes ($n_\text{nod}$) can be of the order of millions. For such large-scale problems direct methods are prohibitive as they do not scale well \cite{Kopanicakova2025}. Efficient alternative methods for large sparse linear systems are those based on Krylov spaces \cite{Saad, Trefethen}, such as conjugate gradients (CG) or the generalised minimal residual (GMRES) methods. However, the convergence of these methods is highly dependent on the condition number of the linear system's matrix \cite{Saad}, which at the same time depends on the domain's discretization. A domain with small and complex geometric features is likely to present large gradients in the discretization, thus resulting in a large condition number \cite{Brenner}. 

A widely used solution to improve the convergence of Krylov methods is to use a preconditioning strategy \cite{Saad, Trefethen}. There are many popular preconditioners, from the relatively simple Jacobi, Gauss-Seidel, incomplete LU (ILU), or successive over-relaxation (SOR) \cite{Meijerink, Saad, Young}, to the more complex methods based on domain-decomposition \cite{Toselli} or multigrid (MG) techniques \cite{Trottenberg}. Relatively simple preconditioners such as ILU or SOR are utilized when there is no geometric information available. On the other hand, domain decomposition methods or MG methods are commonly employed when there is such geometric information available. Of particular interest for linear systems derived from PDEs are the MG methods due to their scalability and favorable computational complexity. MG methods are based on solving the underlying PDEs on a smaller subspace, e.g. a coarser mesh generated from the same geometry, and then use this solution to improve the solution approximation from the finer mesh. \textit{State-of-the-art} MG preconditioners may also lead to significant improvements in convergence, vastly surpassing that of standard preconditioners \cite{xu2017algebraic}. MG methods can be interpreted in the context of subspace correction methods \cite{Xu}, which involve two main algorithmic components: transfer operators and subspace solvers. These are problem-dependent, and substantially affect the applicability and convergence properties of the resulting iterative method \cite{Kopanicakova2025}.

In the last few years, machine learning (ML) methods have been used to solve the underlying PDEs directly, completely bypassing the need for standard numerical methods, such as the aforementioned FEM. Among the most popular ML-based approaches to solve PDEs are physics-informed neural networks (PINNs) \cite{Cai1, Cai2, Raissi} and neural operators (NO) \cite{Jin, Li2021, Lu2021}. PINNs are trained in order to minimize the mean square error of the PDE residual, boundary/initial conditions, and observed data. NOs learn operators, i.e. mappings between function spaces, thus they learn the mapping between the parametric input functions and the corresponding output functions (i.e. the solutions of the PDEs). A substantial advantage of NOs is their reusability; they can be used for solution inference in previously unseen domain locations without extra training efforts. 

A recent strategy is to hybridize classical iterative linear solvers with pretrained NOs by employing hybrid preconditioners \cite{Kahana2023, Lee20251, Lee20252, Zhang2024}. The idea behind this hybridization is to let classical iterative solvers reduce high-frequency modes \cite{Saad} whilst the NOs reduce low-frequency modes due to their spectral bias \cite{Rahaman2019}. These hybrid preconditioners improve the convergence of classical iterative solvers in structured domains compared to when standard preconditioners are used \cite{Lee20251, Lee20253, Kopanicakova2025}. A big disadvantage of these hybrid preconditioning strategies is that they are restricted to the same domain where the corresponding NO was trained, meaning that geometric transferability is still an open issue. Furthermore, even if the convergence is better than that of standard preconditioners in structured domains, it is not clear when complex unstructured domains are used. 

\begin{figure}
\centering
\begin{tikzpicture}[thick, fill opacity=1, text opacity=1]

    \node[rotate=0,minimum width=0.25cm, inner sep=1pt] (training dataset) at (-3,1.5) {\small Training dataset};
    \node[inner sep=0pt] (img1) at (-3.0,-0.35)
    {\includegraphics[width=0.175\textwidth]{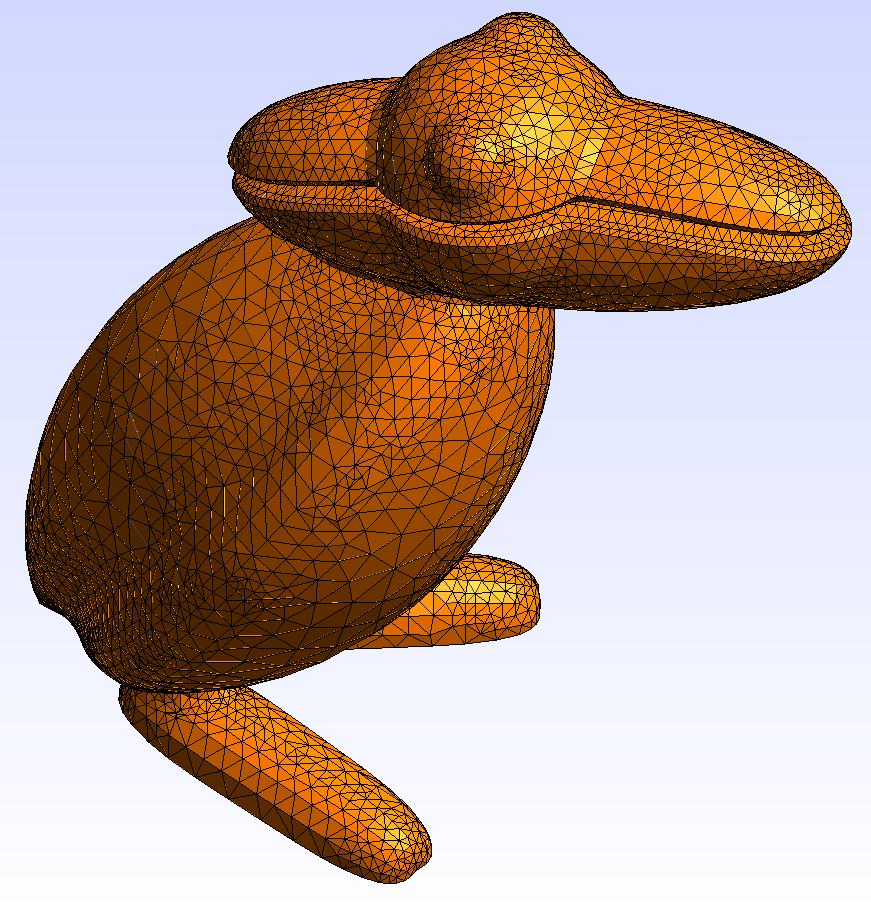}};
    \draw[fill opacity=1.0] (-4.5, -1.85) rectangle (-1.5, 1.2) node[midway]{};

    \draw [line width=2mm, arrows = {-Latex[width'=5pt 0.5, length=15pt]}] (-1.4,-0.35) -- (-0.3,-0.35);

	\draw[fill opacity=0.5, ultra thick] (0, 0) rectangle (0.5, 0.5) node[midway]{$\mathbf{f}$};
    \draw[fill opacity=0.5, ultra thick] (0, -0.6) rectangle (0.5, -0.1) node[midway]{$\mathbf{X}$};
    \draw[fill opacity=0.5, ultra thick] (0, -1.2) rectangle (0.5, -0.7) node[midway]{$\mathbf{d}$};

    \draw[fill opacity=0.5, rounded corners] (1.3, 0) rectangle (2.5, 0.5) node[midway]{\textbf{Branch}};
    \draw[fill opacity=0.5, rounded corners] (3.3, 0) rectangle (4.5, 0.5) node[midway]{\textbf{Scaling}};
    
    \draw[fill opacity=0.5, rounded corners] (1.3, -0.6) rectangle (2.5, -0.1) node[midway]{\textbf{Trunk}};

    \draw[fill opacity=0.5, rounded corners] (1.3, -1.2) rectangle (2.5, -0.7) node[midway]{\textbf{BC}};

    \draw [arrows = {-Latex[width'=0pt .5, length=8pt]}] (0.5,0.25) -- (1.3,0.25);
    \draw [-] (2.5,0.25) -- (3.3,0.25);

    \draw [arrows = {-Latex[width'=0pt .5, length=8pt]}] (0.5,-0.35) -- (1.3,-0.35);

    \draw [arrows = {-Latex[width'=0pt .5, length=8pt]}] (0.5,-0.95) -- (1.3,-0.95);

    \draw [-] (2.5,-0.35) -- (3.5,-0.35);
    \draw [-] (2.5,-0.95) -- (3.5,-0.95);

    \node[rotate=0,minimum width=0.25cm, inner sep=1pt] (odot) at (3.5,-0.65) {\small $\odot$};

    \draw [-] (3.5,-0.35) -- (odot);
    \draw [-] (3.5,-0.95) -- (odot);

    \node[rotate=0,minimum width=0.25cm, inner sep=1pt] (otimes) at (5.25,-0.2) {\small $\otimes $};

    \draw [-] (4.5,0.25) -- (5.25,0.25);
    \draw [-] (5.25,0.25) -- (otimes);

    \draw [-] (odot) -- (5.25,-0.65);
    \draw [-] (5.25,-0.65) -- (otimes);

    \draw[fill opacity=0.5, ultra thick] (6, -0.45) rectangle (6.5, 0.05) node[midway]{$\widetilde{\mathbf{u}}$};

    \draw [arrows = {-Latex[width'=0pt .5, length=8pt]}] (otimes) -- (6,-0.2);

    \draw[fill opacity=1.0, rounded corners] (-0.2, -1.4) rectangle (6.7, 0.7) node[midway]{};
    \node[rotate=0,minimum width=0.25cm, inner sep=1pt] (training dataset) at (3.25,1.1) {\small PTFONet};
    \draw [line width=2mm, arrows = {-Latex[width'=5pt 0.5, length=15pt]}] (3.25,-1.475) -- (3.25,-2.9);

    \draw[fill opacity=0.5, rounded corners] (-4.5, -9) rectangle (6.7, -3) node[midway]{};

    \draw[blue, text=blue, fill opacity=0.5, rounded corners] (-3.8, -4.1) rectangle (0.5, -3.6) node[midway]{$\textbf{Jacobi, GS, SOR, \dots}$};
    \node[inner sep=0pt, text=blue] (prec1) at (-1.65,-3.3)
    {$\text{Standard preconditioner } \mathcal{M}^{-1}_1$};
    \draw [blue, line width=0.5mm, arrows = {-Latex[width'=2pt 0.5, length=8pt]}] (0.5,-3.8) -- (3,-3.8);
    
     \draw[fill opacity=0.5, rounded corners] (3, -4.3) rectangle (6, -3.6) node[midway]{$\left\{ \mathcal{M}^{-1}_1, \mathcal{M}^{-1}_2 \right\}$};
    \node[inner sep=0pt] (prec1) at (4.5,-3.3)
    {$\text{Hybrid preconditioner}$};
    \draw [line width=1mm, arrows = {-Latex[width'=4pt 0.5, length=8pt]}] (4.5,-4.3) -- (4.5,-4.8);

    \draw[red, fill opacity=0.5, rounded corners] (-4.3, -8.5) rectangle (6.5, -5.8) node[midway]{};
    \node[text=red, inner sep=0pt] (prec1) at (1.35,-8.75)
    {$\text{NSPOD preconditioner}$};  

    \draw[fill opacity=0.5, ultra thick] (-3.6, -6.5) rectangle (0.75, -6) node[midway]{$\mathbf{f}_i = \mathbf{f} + \eta \, \mathbf{s}_i, \; \mathbf{s}_i \sim \mathcal{U} \left( \mathbf{f}_\text{min}, \mathbf{f}_\text{max} \right)$};
    \draw[fill opacity=0.5] (-3.6, -7.1) rectangle (-1.6, -6.6) node[midway]{$\left\{ \mathbf{f}_1, \mathbf{f}_2, \dots, \mathbf{f}_m \right\}$};
    \draw[fill opacity=0.5] (-2.1, -7.7) rectangle (-1.6, -7.2) node[midway]{$\mathbf{X}$};
    \draw[fill opacity=0.5] (-2.1, -8.3) rectangle (-1.6, -7.8) node[midway]{$\mathbf{d}$};
    
    \draw[fill opacity=0.5, rounded corners] (-0.8, -7.7) rectangle (0.75, -7.2) node[midway]{\textbf{PTFONet}};
    
    \draw[fill opacity=0.5] (2.2, -6.5) rectangle (5.7, -6) node[midway]{$\mathbf{U} = \begin{bmatrix}
        \widetilde{\mathbf{u}}_1 & \widetilde{\mathbf{u}}_2 & \dots & \widetilde{\mathbf{u}}_m
    \end{bmatrix}$};
    \draw[fill opacity=0.5, rounded corners] (3.5, -7.1) rectangle (4.4, -6.6) node[midway]{\textbf{POD}};
    \draw[fill opacity=0.5] (2.2, -7.7) rectangle (5.7, -7.2) node[midway]{$\mathbf{P}^\text{T} = \begin{bmatrix}
        \pmb{\phi}_1 & \pmb{\phi}_2 & \dots & \pmb{\phi}_m
    \end{bmatrix}$};
    \draw[fill opacity=0.5, ultra thick] (2.2, -8.3) rectangle (5.7, -7.8) node[midway]{$\mathcal{M}^{-1}_2 = \mathbf{P} \left( \mathbf{P}^\text{T} \mathbf{K} \mathbf{P} \right) \mathbf{P}^\text{T}$};

    \draw [-] (-3.6,-6.25) -- (-4.1,-6.25);
    \draw [-] (-4.1,-6.25) -- (-4.1,-6.85);
    \draw [arrows = {-Latex[width'=0pt .5, length=8pt]}] (-4.1,-6.85) -- (-3.6,-6.85);

    \draw [-] (-1.6,-6.85) -- (-1.3,-6.85);
    \draw [-] (-1.6,-8.05) -- (-1.3,-8.05);
    \draw [-] (-1.3,-6.85) -- (-1.3,-8.05);
    \draw [arrows = {-Latex[width'=0pt .5, length=8pt]}] (-1.6,-7.45) -- (-0.8,-7.45);

    \draw [-] (0.75,-7.45) -- (1.05,-7.45);
    \draw [-] (1.05,-7.45) -- (1.05,-6.25);
    \draw [arrows = {-Latex[width'=0pt .5, length=8pt]}] (1.05,-6.25) -- (2.2,-6.25);

    \draw [-] (5.7,-6.25) -- (6.2,-6.25);
    \draw [-] (6.2,-6.25) -- (6.2,-6.85);
    \draw [arrows = {-Latex[width'=0pt .5, length=8pt]}] (6.2,-6.85) -- (4.4,-6.85);

    \draw [-] (3.5,-6.85) -- (1.7,-6.85);
    \draw [-] (1.7,-6.85) -- (1.7,-7.45);
    \draw [arrows = {-Latex[width'=0pt .5, length=8pt]}] (1.7,-7.45) -- (2.2,-7.45);

    \draw [-] (5.7,-7.45) -- (6.2,-7.45);
    \draw [-] (6.2,-7.45) -- (6.2,-8.05);
    \draw [arrows = {-Latex[width'=0pt .5, length=8pt]}] (6.2,-8.05) -- (5.7,-8.05);
    
    \draw [red, -, line width=0.5mm] (1.1,-5.8) -- (1.1,-4.075);
    \draw [red, line width=0.5mm, arrows = {-Latex[width'=2pt 0.5, length=8pt]}] (1.1,-4.1) -- (3,-4.1);

    \draw[green, text=green, fill opacity=0.5, rounded corners] (3.3, -5.6) rectangle (5.7, -4.8) node[midway]{$\textbf{Krylov method}$};

    \node[rotate=0,minimum width=0.25cm, inner sep=1pt] (training dataset) at (-7.5,1.4) {\small Target problem};
    \node[inner sep=0pt] (img1) at (-7.5,-0.35)
    {\includegraphics[width=0.25\textwidth]{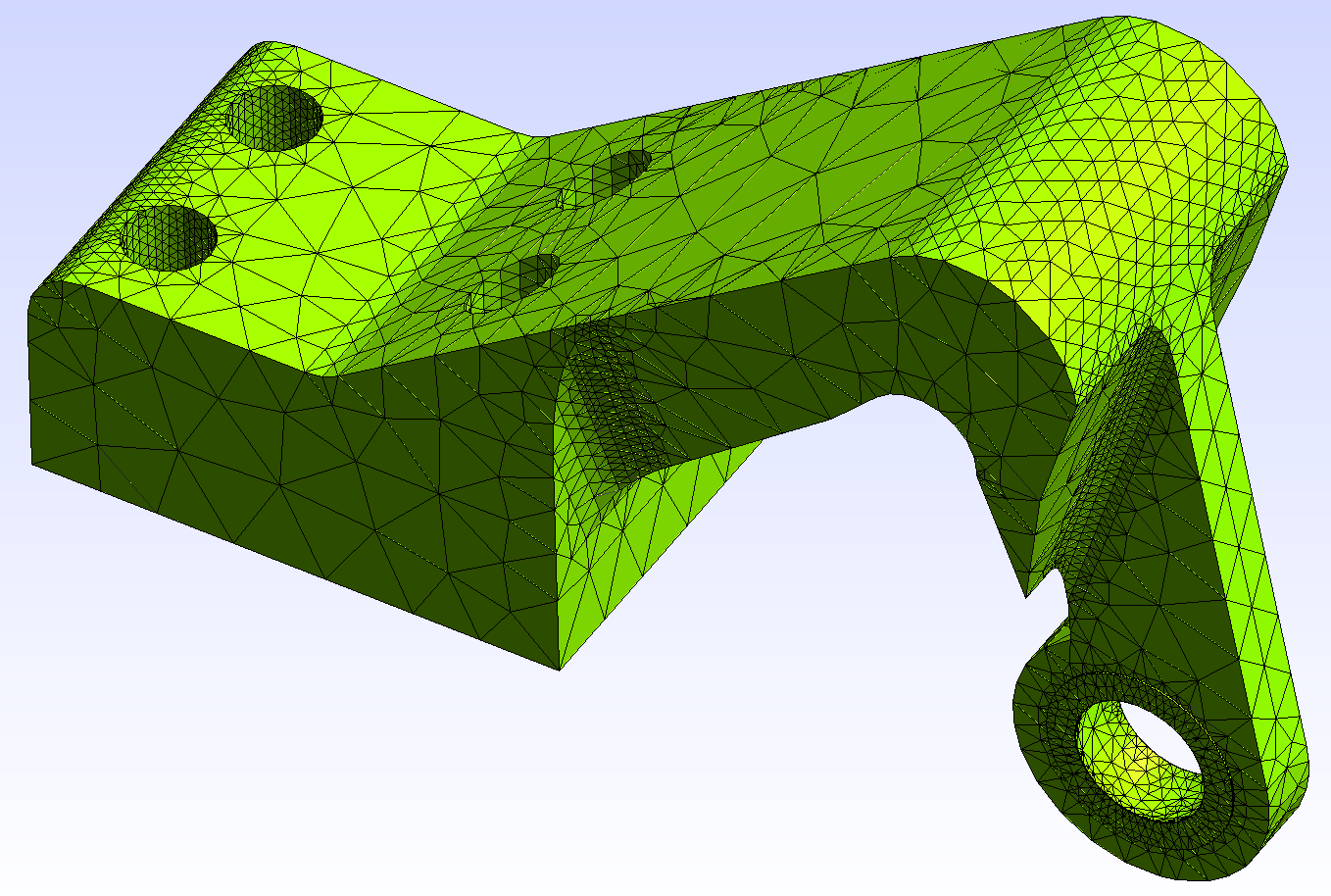}};
    \draw[fill opacity=1.0] (-9.6, -1.75) rectangle (-5.4, 1.1) node[midway]{};
    \node[rotate=0,minimum width=0.25cm, inner sep=1pt] (training dataset) at (-7.5,-2.6) {\small $\begin{cases}
        \nabla \cdot \pmb{\sigma} + \mathbf{b} = \mathbf{0}, &\text{ in } \Omega_0, \\
        \pmb{\sigma} \mathbf{n} - \mathbf{t} = \mathbf{0}, &\text{ on } \partial\Omega_0^N, \\
        \mathbf{u} - \overline{\mathbf{u}} = \mathbf{0}, &\text{ on } \partial\Omega_0^D,
    \end{cases}$};

    \draw [line width=2mm, arrows = {-Latex[width'=5pt 0.5, length=15pt]}] (-7.5,-3.5) -- (-7.5,-5.55);
    
    \draw [line width=2mm, arrows = {-Latex[width'=5pt 0.5, length=15pt]}] (-4.6,-6) -- (-6.2,-6);
    \draw[fill opacity=0.5] (-8.6, -5.65) rectangle (-6.4, -6.35) node[midway]{\LARGE $\mathbf{K} \mathbf{u} = \mathbf{f}$};
    \draw[fill opacity=0.5] (-8, -8) rectangle (-7, -8.7) node[midway]{\LARGE $\mathbf{u}^\star$};
    \draw [line width=2mm, arrows = {-Latex[width'=5pt 0.5, length=15pt]}] (-7.5,-6.425) -- (-7.5,-7.9);

	\end{tikzpicture}
	\caption{A graphical overview of the NSPOD hybrid preconditioned Krylov-based linear iterative solvers. The PTFONet is trained on CAD geometries in an \textit{offline} stage, encoding the geometries used during training, which allows for generation of preconditioners for the same geometry but different input functions, or for completely different geometries altogether, in the subsequent \textit{online} stage. This hybrid preconditioner is applied to a new combination of input functions in an \textit{online} stage in which the preconditioner is created and applied.}
	\label{fig:overview}
\end{figure}
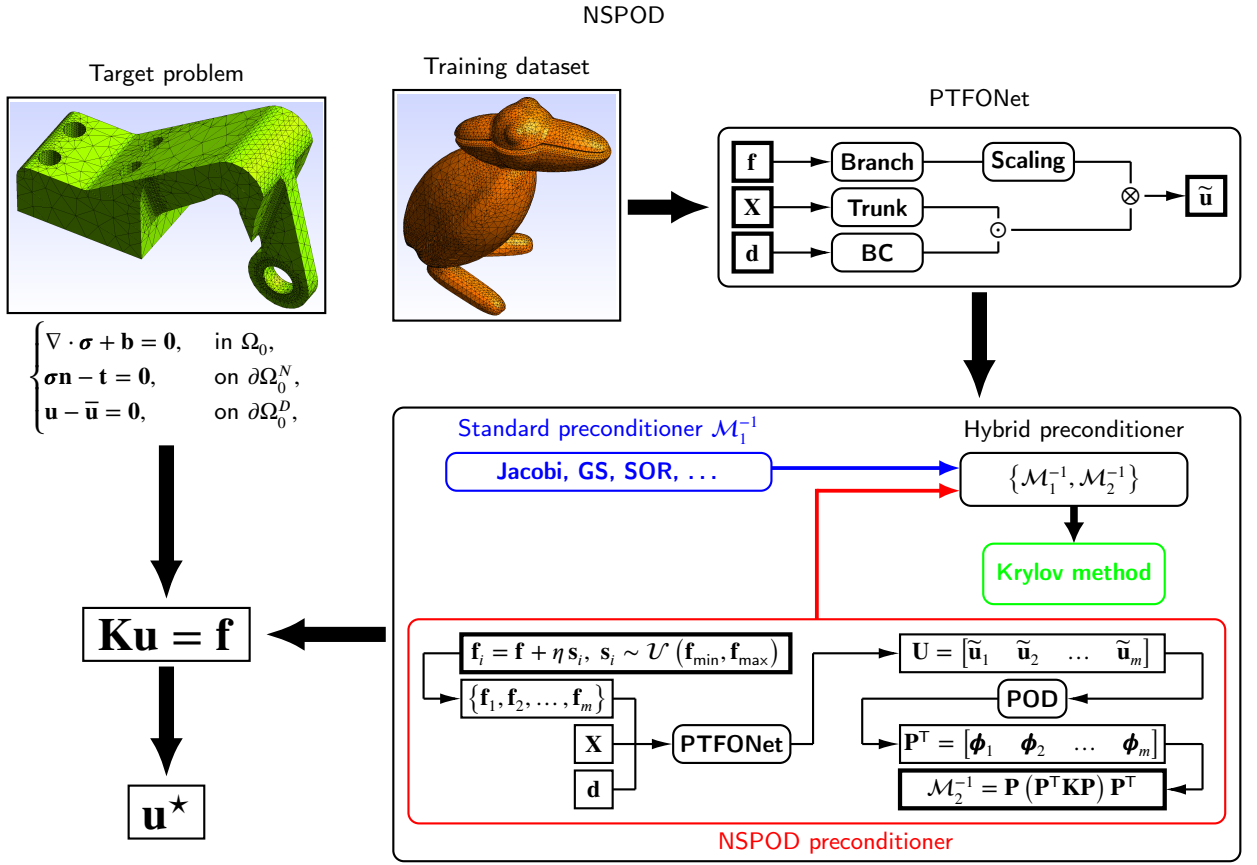

In this study we propose hybrid iterative solvers with neural preconditioners to solve parametric PDEs on unstructured domains. We use deep operator networks (DeepONet) as our NOs, as it performs better on unstructured domains compared to other NO architectures \cite{Lu2021, Lu2022}. We introduce PointTransformer \cite{PointTransformer} layers in modified architectures for both branch and trunk subnetworks of DeepONet: the branch network is based on the encoder portion of the semantic segmentation variant of PointTransformer networks, whilst the trunk network consists of several stacked layers of PointTransformer blocks \cite{PointTransformer}. PointTransformers are powerful architectures that adapt traditional transformers to point cloud data; they are able to capture long-range dependencies between points in 3D space. Finally, at the end of the branch subnetwork we also include a customized squeezed-and-excitation network, producing a scaled output acting as an attention mechanism \cite{Hu2018, Lee20251, Vaswani2023}. By combining all these components we build PTFONet. We train these PTFONets on a number of CAD geometries, several of them if we want to achieve a certain degree of geometric transferability. We then use these trained networks to perform solution inference and build preconditioners based on Proper Orthogonal Decomposition (POD) of these inferred solutions; we denote this hybrid preconditioner as Neural Subspace Proper Orthogonal Decomposition (NSPOD). A sketch of this methodology is shown in Fig. \ref{fig:overview}. The superior performance of the resulting preconditioners and the assessment of their geometric transferability is shown via various numerical experiments.

This paper is organized as follows: in Section \ref{sec:model}, we describe the studied model problem and review the FEM applied to elliptic systems; in Section \ref{sec:linear_solvers}, we introduce the hybrid preconditioning framework that combines the pretrained PTFONet with Krylov-based linear iterative solvers; in Section \ref{sec:ptfonet} the mathematical framework for the newly developed version of DeepONet is thoroughly described; in Sections \ref{sec:network_architecture} and \ref{sec:training} we describe the employed network architecture framework and how the networks are trained, respectively; in Section \ref{sec:neural_POD} we describe our newly develop hybrid preconditioner NSPOD; in Section \ref{sec:computational_cost} we describe the computational cost of the developed hybrid preconditioner for both \textit{offline} and \textit{online} stages. Finally, in Section \ref{sec:results}, we demonstrate the numerical performance of the proposed hybrid preconditioner strategy via numerical experiments. Discussion, conclusion, limitations, and future work are all discussed in Section \ref{sec:conclusion}.

%%%%%%%%%%%%%%%%%%%%%%%%%%%%%%%%%%%%%%%%%%%

\section{Model problem} \label{sec:model}
\noindent In this manuscript we focus on the study of linearized equations in quasi-static solid mechanics. A body occupies a bounded region $\Omega_0 \subset \mathbb{R}^{n_\text{dim}}$, with boundary $\partial\Omega_0$ which is not necessarily simply connected nor convex; $n_\text{dim}$ is the number of spatial dimensions considered. This body is subjected to surface tractions, or \textit{Neumann boundary conditions} (BCs), on $\partial\Omega_0^N$, and to prescribed displacements, or \textit{Dirichlet boundary conditions}, on $\partial\Omega_0^D$, such that $\partial\Omega_0 = \partial\Omega_0^N \cup \partial\Omega_0^D$ and $\partial\Omega_0^N \cap \partial\Omega_0^D = \emptyset$.

By considering a quasi-static regime, the mechanical behavior of the considered body is characterized by the \textit{balance of linear momentum} equation, such that
\begin{align}
\nabla \cdot \pmb{\sigma} + \mathbf{b} = \mathbf{0}, &\text{ in } \Omega_0, \label{eq:balance_linear_momentum} \\
\pmb{\sigma} \mathbf{n} - \mathbf{t} = \mathbf{0}, &\text{ on } \partial\Omega_0^N, \label{eq:neumann_bcs} \\
\mathbf{u} - \overline{\mathbf{u}} = \mathbf{0}, &\text{ on } \partial\Omega_0^D, \label{eq:dirichet_bcs}
\end{align}

\noindent where $\mathbf{u}$ is the displacement field, $\nabla \cdot \left( \cdot \right)$ is the divergence of $\left( \cdot \right)$, $\pmb{\sigma} = \mathbb{D}^e : \pmb{\varepsilon}$ is the Cauchy's stress tensor, $\mathbb{D}^e$ is the fourth-order constant elasticity tensor, $\pmb{\varepsilon}$ is the infinitesimal strain tensor, $\mathbf{b}$ is the volumetric body force, $\mathbf{n}$ is the outward normal defined on $\partial\Omega_0^N$, $\mathbf{t}$ is the surface traction, and $\overline{\mathbf{u}}$ is the prescribed displacement field on $\partial\Omega_0^D$.

As a preliminary study, and in order to simplify learning, we assume the following here: homogeneous Dirichet BCs, which are applied on a strict non-empty subset of $\partial\Omega_0$, namely $\partial\Omega_0^D \subset \partial\Omega_0$, such that $0 < \frac{\left| \partial\Omega_0^D \right|}{\left| \partial\Omega_0 \right|} \leq p < 1 \left( \text{thus } \partial\Omega_0^D \neq \emptyset \text{ and } \partial\Omega_0^D \neq \partial\Omega_0 \right)$; homogeneous Neumann BCs which are applied to the whole $\partial\Omega_0^N$; isotropic elasticity with a spatially-constant Young's modulus $\left( E \right)$, a fixed Poisson's ratio $\left( \nu \right)$ set to a predetermined value, and spatially-constant body forces $\left( \mathbf{b} \right)$. In the following we also assume that the well-posedness of the considered system is ensured. Now let $\mathcal{V} \subset L^2 \left( \Omega_0 \right)$ be a specific Hilbert space consisting of functions on $\Omega_0$. Therefore, Eqs. \ref{eq:balance_linear_momentum}, \ref{eq:neumann_bcs} and \ref{eq:dirichet_bcs} admit the following weak formulation: Find $\mathbf{u} \in \mathcal{V}$ such that
\begin{equation}
\int_{\Omega_0} \left( \mathbb{D}^e : \pmb{\varepsilon} \right) : \nabla\pmb{\eta} \, \text{d}V - \int_{\Omega_0} \mathbf{b} \cdot \pmb{\eta} \, \text{d}V = 0, \forall \pmb{\eta} \in \mathcal{V}. \label{eq:weak_form}
\end{equation}

\noindent Now consider the parameter space $\mathcal{H} = \left\{ E, p, \mathbf{b}, \sigma \right\} \subset \mathbb{R}^{2 + n_\text{dim}} \times S_k$, where $\sigma \in S_k$ is a permutation of elements $n = \left\{ n_1, n_2, \dots, n_k \right\} \in \mathbb{N}_k$; the meaning of this permutation will be revealed when the dataset creation process is described in Section \ref{sec:training}. By considering this parameter space, Eq. \ref{eq:weak_form} can be written as: Given $\pmb{\xi} \in \mathcal{H}$, find $\mathbf{u} \in \mathcal{V}$ such that
\begin{equation} \label{eq:weak_parameterized_form}
\mathcal{A} \left( \mathbf{u}, \pmb{\eta}; \pmb{\xi} \right) = \left( \mathbf{b}, \pmb{\eta}; \pmb{\xi} \right), \forall \pmb{\eta} \in \mathcal{V},
\end{equation}

\noindent where $\mathcal{A} \left( \cdot, \cdot; \pmb{\xi} \right)$ and $\left( \cdot, \cdot; \pmb{\xi} \right)$ denote a bilinear form obtained from $\nabla \cdot \pmb{\sigma}$ and an inner product defined on $L^2 \left( \Omega_0 \right)$, respectively. Equation \ref{eq:weak_parameterized_form} can be solved via the FEM, yielding a system of linear equations
\begin{equation} \label{eq:linear_system}
\mathbf{K} \mathbf{u} = \mathbf{f},
\end{equation}

\noindent where $\mathbf{u}$ are the nodal displacements $\mathbf{u} \in \mathbb{R}^{n_\text{nod} \times n_\text{dim}}$, $\mathbf{K} \in \mathbb{R}^{\left( n_\text{nod} \times n_\text{dim} \right) \times \left( n_\text{nod} \times n_\text{dim} \right)}$ is the stiffness matrix, and $\mathbf{f} \in \mathbb{R}^{n_\text{nod} \times n_\text{dim}}$ is the load vector. In here both $\mathbf{K}$ and $\mathbf{f}$ depend on the parameters $\pmb{\xi}$, and their components are computed as
\begin{equation*}
\mathbf{K}_{ij} = \mathcal{A} \left( \phi_i, \phi_j; \pmb{\xi} \right), \, \mathbf{f}_i = \left( \mathbf{b}, \phi_i; \pmb{\xi} \right), 1 \leq i, j \leq n_\text{nod} \times n_\text{dim},
\end{equation*}

\noindent where $\left\{ \phi_i \right\}_{i = 1}^{n_\text{nod} \times n_\text{dim}}$ denotes a set of nodal basis functions on a suitable finite element space $\mathcal{V}_h$ defined on a quasi-uniform tetrahedral mesh $\mathcal{T}^h$ of $\Omega_0$, $h$ being the mesh size. To guarantee convergence of the solution to a continuous solution, coercivity and elliptic regularity need to be assumed \cite{Brenner}.

%%%%%%%%%%%%%%%%%%%%%%%%%%%%%%%%%%%%%%%%%%%

\section{Methodology} 
\noindent In order to solve the resulting system of linear Equations (Eq. \ref{eq:linear_system}), various preconditioned iterative solvers might be employed, with Krylov-based solvers being the most widely used~\cite{Briggs,Saad,Toselli}, especially when the linear system is large. Defining the preconditioner is usually a challenging task. Hybrid preconditioners \cite{Kopanicakova2025, Zhang2024} may assist in simplifying this task by employing approximated inferred solutions in different ways, such as by constructing the bases of projection operators in a MG fashion \cite{Kopanicakova2025, Lee20251}, or by using direct inference when relaxation methods are employed \cite{Kahana2023, Zhang2024}. Despite these efforts, performance on arbitrary, unstructured, domains is still an unresolved issue. In this section we briefly introduce the basic concepts behind hybrid preconditioned Krylov-based linear iterative solvers, and later we propose our novel NSPOD approach, which significantly improves the convergence behavior of even \text{state-of-the-art} MG approaches.

%%%%%%%%%%%%%%%%%%%%%%%%%%%%%%%%%%%%%%%%%%%

\subsection{Hybrid preconditioners} \label{sec:linear_solvers}
\noindent Direct inference through an NO is most commonly employed when relaxation methods are used because Krylov-based methods require linear preconditioners, and inference employing neural networks is generally not linear. However, flexible Krylov-based methods \cite{Simoncini} can be used with nonlinear preconditioners, so we introduce here hybrid preconditioners for relaxation methods using DeepONet. A standard relaxation method can be formulated as a preconditioned Richardson iteration \cite{Saad2}. The $i-$th iteration $\mathbf{u}^{\left( i \right)}$ of the solution given by a standard relaxation method is then
\begin{equation*}
\begin{cases}
\mathbf{r}^{\left( i \right)} &= \mathbf{f} - \mathbf{K} \mathbf{u}^{\left( i \right)}, \\
\mathbf{u}^{\left( i + 1 \right)} &= \mathbf{u}^{\left( i \right)} + \mathcal{M}^{-1} \mathbf{r}^{\left( i \right)},
\end{cases}
\end{equation*}

\noindent where $\mathcal{M}^{-1}$ denotes the specific preconditioner of the Richardson iteration, e.g. choosing $\mathcal{M}^{-1}$ as the diagonal or lower triangular part of $\mathbf{K}$ leads, respectively, to the Jacobi or Gauss-Seidel methods. The hybrid preconditioned Richardson iteration, also known as HINTS \cite{Kahana2023, Zhang2024}, can be written as
\begin{equation}
\begin{cases}
\mathbf{r}^{\left( i \right)} &= \mathbf{f} - \mathbf{K} \mathbf{u}^{\left( i \right)}, \\
\mathbf{u}^{\left( i + \frac{1}{2} \right)} &= \mathbf{u}^{\left( i \right)} + \mathcal{M}_1^{-1} \left( \mathbf{r}^{\left( i \right)} \right), \\
\mathbf{r}^{\left( i + \frac{1}{2} \right)} &= \mathbf{f} - \mathbf{K} \mathbf{u}^{\left( i + \frac{1}{2} \right)}, \\
\mathbf{u}^{\left( i + 1 \right)} &= \mathbf{u}^{\left( i + \frac{1}{2} \right)} + \mathcal{M}_2^{-1} \left( \mathbf{r}^{\left( i + \frac{1}{2} \right)} \right),
\end{cases}
\end{equation}

\noindent where $\mathcal{M}_1^{-1}$ and $\mathcal{M}_2^{-1}$ denote, respectively, the $n_r$ iterations of the relaxation method and the prediction of the pretrained NO (DeepONets are commonly used here). More details about this form of hybrid preconditioner may be found in \cite{Lee20251, Zhang2024}.

Hybrid preconditioning of standard Krylov-based methods can be achieved by the trunk basis (TB) approach \cite{Kopanicakova2025}. This method extracts the basis functions from the last layer of the trunk network of DeepONet and uses them to construct a linear preconditioner. We introduce here the TB-based prolongation operator $\mathbf{P} \in \mathbb{R}^{n_\text{nod} \times n_{bf}}$, where $n_{bf}$ is the number of basis functions used. The $\left(  i, j \right)-$th component of $\mathbf{P}$ is $\mathbf{P}_{ij} = \mathcal{T}_j \left( \mathbf{X}_i \right)$, where $\mathcal{T}_j \left( \mathbf{X}_i \right)$ is the $j-$th component of the last layer of the trunk network evaluated at the spatial location $\mathbf{X}_i \in \Omega_0 \subset \mathbb{R}^3$. To define a well-conditioned coarse-level operator $\mathbf{K}_c = \mathbf{P}^\text{T} \mathbf{K} \mathbf{P}$, it is required that $\mathbf{P}$ has full rank and orthogonal columns, which is straightforwardly achieved by performing a QR decomposition of $\mathbf{P}$ and using $\mathbf{Q}$ as the prolongation of operator instead \cite{Lee20251}. Therefore, the TB-based preconditioner $\mathcal{M}_\text{TB}^{-1}$ is defined as
\begin{equation} \label{eq:krylov_hybrid_preconditioner}
\mathcal{M}_\text{TB}^{-1} \coloneqq \mathbf{P} \left( \mathbf{P}^\text{T} \mathbf{K} \mathbf{P} \right) \mathbf{P}^\text{T} = \mathbf{P} \mathbf{K}_c \mathbf{P}^\text{T},
\end{equation}

\noindent where $\mathbf{K}_c$ is the coarse-level problem commonly used in MG methods \cite{Briggs, Trottenberg}. It is clear that the linear, symmetric, operator in Eq. \ref{eq:krylov_hybrid_preconditioner} can be used as the preconditioner in standard Krylov-based methods. Adding everything together, given the residual $\mathbf{r}^{\left( i \right)} = \mathbf{f} - \mathbf{K} \mathbf{u}^{\left( i \right)}$ in a Krylov iteration, the hybrid preconditioned residual is defined as
\begin{equation}
\begin{cases}
\mathbf{u}' &= \mathbf{u}^{\left( i \right)} + \mathcal{M}_1^{-1} \left( \mathbf{r}^{\left( i \right)} \right), \\
\mathbf{u}'' &= \mathbf{u}' + \mathcal{M}_\text{TB}^{-1} \left( \mathbf{f} - \mathbf{K} \mathbf{u}' \right), \\
\mathbf{u}''' &= \mathbf{u}'' + \mathcal{M}_1^{-1} \left( \mathbf{f} - \mathbf{K} \mathbf{u}'' \right), \\
\mathbf{z}^{\left( i \right)} &= \mathbf{f} - \mathbf{K} \mathbf{u}''',
\end{cases}
\end{equation}

\noindent where $\mathcal{M}_1^{-1}$ denotes the $n_r$ iterations of the relaxation method used in this MG-like two-level preconditioner. More details of this preconditioner, as well as a schematic representation, can be seen in \cite{Lee20251}.

%%%%%%%%%%%%%%%%%%%%%%%%%%%%%%%%%%%%%%%%%%%

\subsection{PointTransformer DeepONet} \label{sec:ptfonet}
\noindent We first introduce the ``vanilla'' DeepONet architecture here and briefly describe its mathematical formulation; a more detailed description can be found in \cite{Lee20251}. We then introduce our newer architecture, PointTransformer DeepONet (PTFONet).

Let $\mathcal{X}$ and $\mathcal{Y}$ be infinite dimensional Banach spaces, and let $\mathcal{G} : \mathcal{X} \rightarrow \mathcal{Y}$ be some operator to be approximated via DeepONet, with the corresponding approximation being $\mathcal{G}_\text{NN}$. The output of DeepONet is defined as a contraction between the outputs of two different neural networks, called the branch and the trunk networks. The branch network $\mathcal{B} : \mathbb{R}^m \rightarrow \mathbb{R}^p$ is a vector-valued network of dimension $p$ that takes a discretized function $\mathbf{f}_m \in \mathcal{X}_m$ as the input. We assume that $\mathcal{X}_m$ is an $m-$dimensional space satisfying
\begin{equation*}
\forall f \in \mathcal{X}, \exists \mathbf{f}_m \in \mathcal{X}_m \text{ such that } \lim_{m \to \infty} \mathbf{f}_m = f.
\end{equation*}

\noindent The trunk network $\mathcal{T} : \Omega_0 \subset \mathbb{R}^{n_\text{dim}} \rightarrow \mathbb{R}^p$ is another vector-valued network of dimension $p$, but it now takes coordinates in the domain $\Omega_0$ as the input. The final output of DeepONet is then given by
\begin{equation*}
\mathcal{G}_\text{NN} \left( \mathbf{f} \right) \left( \mathbf{X} \right) \coloneqq \langle \mathcal{B} \left( \mathbf{f} \right), \mathcal{T} \left( \mathbf{X} \right) \rangle, \; \forall \mathbf{f} \in \mathcal{X}_m, \mathbf{X} \in \Omega_0,
\end{equation*}

\noindent where $\mathbf{X}$ are spatial coordinates and $\langle \cdot_1, \cdot_2 \rangle$ is the inner product between $\cdot_1$ and $\cdot_2$. With this setting, the universal approximation theorem for operators \cite{Lu2021} guarantees that DeepONet may sufficiently approximate any target operator in $L^2-$norm.

Given the quasi-uniform tetrahedral discretization $\mathcal{T}^h$, the adjacency matrix $\mathbf{A}$ can be suitably defined, which consists of $\left( i, j \right)$ components equal to one if the $i-$th and $j-$th nodes are connected and zero otherwise. Alternatively, one can use the set of edges $\mathcal{E} \subseteq \left\{ \left( v_1, v_2 \right)_1, \left( v_1, v_2 \right)_2, \dots | v_1, v_2 \in \left\{ 1, 2, \dots, n_\text{nod} \right\} \right\} \eqqcolon \mathcal{N}$. Our new architecture, PTFONet, utilizes PointTransformers \cite{PointTransformer} as the main convolution layers in both branch and trunk networks, and it requires of spatial coordinates and connectivity information such as that provided by the coordinates $\mathbf{X}$ and edge indices $\mathcal{E}$. PTFONet also adds an additional subnetwork, BC network from now onwards, that deals with the position of fully restrained nodes in the mesh, which correspond to the locations where homogeneous Dirichlet BCs have been applied. We denote the new versions of the branch, trunk and BC networks as $\widetilde{\mathcal{B}}$, $\widetilde{\mathcal{T}}$, and $\widetilde{\mathcal{D}}$, respectively. In addition, a scaling network $\mathcal{S}$, inspired by the squeeze-and-excitation network \cite{Hu2018}, is immediately added after the new branch network in order to enhance the performance of PTFONet; more information about this scaling network can be found in \cite{Lee20251}. 

Therefore PTFONet consists of four subnetworks: $\widetilde{\mathcal{B}} : \mathbb{R}^m \times \mathbb{R}^{n_\text{dim}} \times \mathcal{N} \rightarrow \mathbb{R}^{n_\text{dim} \times p}$, $\widetilde{\mathcal{T}} : \mathbb{R}^{n_\text{dim}} \times \mathcal{N} \rightarrow \mathbb{R}^p$, $\widetilde{\mathcal{D}} : \left\{ 0, 1 \right\} \times \mathbb{R}^{n_\text{dim}} \times \mathcal{N} \rightarrow \mathbb{R}^p$, and $\mathcal{S} : \mathbb{R}^{n_\text{dim} \times p} \rightarrow \mathbb{R}^{n_\text{dim} \times p}$. First the output of the branch network gets passed to the scaling network, yielding $\mathcal{S} \circ \widetilde{\mathcal{B}} : \mathbb{R}^m \times \mathbb{R}^{n_\text{dim}} \times \mathcal{N} \rightarrow \mathbb{R}^{n_\text{dim} \times p}$, where $\circ$ is a function composition; then, the outputs of the trunk and BC network get element-wise multiplied, yielding $\widetilde{\mathcal{T}} \odot \widetilde{\mathcal{D}} : \mathbb{R}^{n_\text{dim}} \times \left\{ 0, 1 \right\} \times \mathcal{N} \rightarrow \mathbb{R}^p$, where $\odot$ is the Hadamard product; finally, these two latter outputs get contracted appropriately such that $\left( \mathcal{S} \circ \widetilde{\mathcal{B}} \right) \otimes \left( \widetilde{\mathcal{T}} \odot \widetilde{\mathcal{D}} \right) : \mathbb{R}^m \times \mathbb{R}^{n_\text{dim}} \times \left\{ 0, 1 \right\} \times \mathcal{N} \rightarrow \mathbb{R}^{n_\text{dim}}$, where $\otimes$ here is the contraction that yields an output of the desired dimensions. The output of PTFONet, namely $\mathcal{G}_\text{PTF}$, is then given by
\begin{equation} \label{eq:ptfo_inference}
\mathcal{G}_\text{PTF} \left( \mathbf{f}, \mathbf{X}, \mathcal{E} \right) \left( \mathbf{X}, d, \mathcal{E} \right) \coloneqq \langle \mathcal{S} \circ \widetilde{\mathcal{B}} \left( \mathbf{f}, \mathbf{X}, \mathcal{E} \right), \widetilde{\mathcal{T}} \odot \widetilde{\mathcal{D}} \left( \mathbf{X}, d, \mathcal{E} \right) \rangle, \; \forall \mathbf{f} \in \mathcal{X}_m, \; \mathbf{X} \in \Omega_0, \; d \in \left\{ 0, 1 \right\}, \; \mathcal{E} \subseteq \mathcal{N}.
\end{equation}

\noindent In the implementation of $\widetilde{\mathcal{D}}$, the position of homogeneous Dirichlet BCs $\overline{\mathbf{X}} \in \partial\Omega^D_0$ gets encoded as a binary variable $d$ which is defined as for all $\mathbf{X} \in \Omega_0$,
\begin{equation*}
d = 
\begin{cases}
1 &\text{ if } \mathbf{X} \in \partial\Omega^D_0, \\
0 &\text{ if } \mathbf{X} \notin \partial\Omega^D_0.
\end{cases}
\end{equation*}

\noindent It is also relevant to mention that, because we are dealing with linearized solid mechanics, isotropic elasticity, and both Dirichlet and Neumann BCs are homogeneous, the relevant input parameters to this study are the Young's modulus $E$ and the volumetric forces $\mathbf{b}$, both considered to be spatially constant. The Poisson's ratio $\nu$ is fixed and assumed to be $\nu = 0.2$.
Thus, $\mathbf{f} = \left\{ E, \mathbf{b} \right\}$ in Eq. \ref{eq:ptfo_inference}, leading to
\begin{equation*}
\mathcal{G}_\text{PTF} \left( E, \mathbf{b}, \mathbf{X}, \mathcal{E} \right) \left( \mathbf{X}, d, \mathcal{E} \right) \coloneqq \langle \mathcal{S} \circ \widetilde{\mathcal{B}} \left( E, \mathbf{b}, \mathbf{X}, \mathcal{E} \right), \widetilde{\mathcal{T}} \odot \widetilde{\mathcal{D}} \left( \mathbf{X}, d, \mathcal{E} \right) \rangle, \; \forall E \in \mathbb{R}^+, \; \forall \mathbf{f} \in \mathcal{X}_3, \; \mathbf{X} \in \Omega_0, \; d \in \left\{ 0, 1 \right\}, \; \mathcal{E} \subseteq \mathcal{N}.
\end{equation*}

%%%%%%%%%%%%%%%%%%%%%%%%%%%%%%%%%%%%%%%%%%%

\subsection{Network architecture} \label{sec:network_architecture}
\begin{figure}
\centering
\begin{tikzpicture}

	\node[rotate=90,minimum width=2.75cm] (x1) at (1.25,0) {\small$\left[ \mathbf{E}^h, \mathbf{b}^h \right]$};
	\node[h,rotate=90,minimum width=2.75cm] (mlp1) at (2.75,0) {\small MLP};
	\node[conv,rotate=90,minimum width=2.75cm] (pntb1) at (3.25,0) {\small PointTransformer};
	\node[fc,rotate=90,minimum width=2.75cm] (td1) at (4.75,0) {\small TransformerDown};
	\node[conv,rotate=90,minimum width=2.75cm] (pntb2) at (5.25,0) {\small PointTransformer};
	\node[fc,rotate=90,minimum width=2.75cm] (td2) at (6.75,0) {\small TransformerDown};
	\node[conv,rotate=90,minimum width=2.75cm] (pntb3) at (7.25,0) {\small PointTransformer};
	\node[fc,rotate=90,minimum width=2.75cm] (td3) at (8.75,0) {\small TransformerDown};
	\node[conv,rotate=90,minimum width=2.75cm] (pntb4) at (9.25,0) {\small PointTransformer};
	\node[pool,rotate=90,minimum width=2.75cm] (pool) at (10.75,0) {\small GlobalAvgPooling};
    \node[h,rotate=90,minimum width=2.75cm] (mlp2) at (12.25,0) {\small MLP};
    \node[bn,rotate=90,minimum width=2.75cm] (SE) at (12.75,0) {\small SE};

    \node[view,rotate=0,minimum width=1cm] (b1) at (1.25,-2) {\small $n_\text{nod} \times 4$};
    \node[view,rotate=0,minimum width=1cm] (b2) at (3,-2) {\small $n_\text{nod} \times C_{b1}$};
    \node[view,rotate=0,minimum width=1cm] (b4) at (5,-2) {\small $\frac{n_\text{nod}}{4} \times C_{b2}$};
    \node[view,rotate=0,minimum width=1cm] (b5) at (7,-2) {\small $\frac{n_\text{nod}}{16} \times C_{b3}$};
    \node[view,rotate=0,minimum width=1cm] (b6) at (9,-2) {\small $\frac{n_\text{nod}}{64} \times C_{b4}$};
    \node[view,rotate=0,minimum width=1cm] (b7) at (10.75,-2) {\small $1 \times C_{b4}$};
    \node[view,rotate=0,minimum width=1cm] (b8) at (12.5,-2) {\small $1 \times 3 n_{bf}$};

    \node[rotate=90,minimum width=2.75cm] (x2) at (1.25,-4.5) {\small$\mathbf{X}^h$};
    \node[h,rotate=90,minimum width=2.75cm] (mlp3) at (3,-4.5) {\small MLP};
	\node[conv,rotate=90,minimum width=2.75cm] (pntb5) at (4.75,-4.5) {\small PointTransformer};
    \node[conv,rotate=90,minimum width=2.75cm] (pntb6) at (6.5,-4.5) {\small PointTransformer};
    \node[conv,rotate=90,minimum width=2.75cm] (pntb7) at (8,-4.5) {\small PointTransformer};
    \node[h,rotate=90,minimum width=2.75cm] (mlp4) at (8.5,-4.5) {\small MLP};

    \node[view,rotate=0,minimum width=1cm] (t1) at (1.25,-6.5) {\small $n_\text{nod} \times 3$};
    \node[view,rotate=0,minimum width=1cm] (t2) at (3,-6.5) {\small $n_\text{nod} \times C_{t1}$};
    \node[view,rotate=0,minimum width=1cm] (t3) at (4.75,-6.5) {\small $n_\text{nod} \times C_{t2}$};
    \node[view,rotate=0,minimum width=1cm] (t4) at (6.5,-6.5) {\small $n_\text{nod} \times C_{t3}$};
    \node[view,rotate=0,minimum width=1cm] (t5) at (8.25,-6.5) {\small $n_\text{nod} \times n_{bf}$};

    \node[rotate=90,minimum width=2.75cm] (x3) at (1.25,-9) {\small$\mathbf{d}^h$};
    \node[h,rotate=90,minimum width=2.75cm] (mlp5) at (3,-9) {\small MLP};
	\node[conv,rotate=90,minimum width=2.75cm] (pntb8) at (4.75,-9) {\small PointTransformer};
    \node[conv,rotate=90,minimum width=2.75cm] (pntb9) at (6.5,-9) {\small PointTransformer};
    \node[conv,rotate=90,minimum width=2.75cm] (pntb10) at (8,-9) {\small PointTransformer};
    \node[h,rotate=90,minimum width=2.75cm] (mlp6) at (8.5,-9) {\small MLP};

    \node[view,rotate=0,minimum width=1cm] (t1) at (1.25,-11) {\small $n_\text{nod} \times 1$};
    \node[view,rotate=0,minimum width=1cm] (t2) at (3,-11) {\small $n_\text{nod} \times C_{t1}$};
    \node[view,rotate=0,minimum width=1cm] (t3) at (4.75,-11) {\small $n_\text{nod} \times C_{t2}$};
    \node[view,rotate=0,minimum width=1cm] (t4) at (6.5,-11) {\small $n_\text{nod} \times C_{t3}$};
    \node[view,rotate=0,minimum width=1cm] (t5) at (8.25,-11) {\small $n_\text{nod} \times n_{bf}$};
		
	\draw[->] (x1) -- (mlp1);
	\draw[->] (pntb1) -- (td1);
	\draw[->] (pntb2) -- (td2);
	\draw[->] (pntb3) -- (td3);
    \draw[->] (pntb4) -- (pool);
    \draw[->] (pool) -- (mlp2);
	
	\draw[->] (x2) -- (mlp3);
	\draw[->] (mlp3) -- (pntb5);
	\draw[->] (pntb5) -- (pntb6);
	\draw[->] (pntb6) -- (pntb7);
	
	\draw[->] (x3) -- (mlp5);
	\draw[->] (mlp5) -- (pntb8);
	\draw[->] (pntb8) -- (pntb9);
	\draw[->] (pntb9) -- (pntb10);

    \node[rotate=90,minimum width=0.5cm] (odot) at (10,-6.75) {\Large $\odot$};
    \node[rotate=90,minimum width=0.5cm] (otimes) at (11.5,-6.75) {\Large $\otimes$};
        
	\draw[-] (SE) -- (13.75,0);
	\draw[-] (13.75,0) -- (13.75,-3.25);
	
	\draw[-] (mlp4) -- (10,-4.5);
    \draw[-] (mlp6) -- (10,-9);
	\draw[-] (10,-4.5) -- (odot);
    \draw[-] (odot) -- (10,-9);
    \draw[-] (odot) -- (otimes);
	\draw[-] (13.75,-3.25) -- (11.5,-3.25);
	\draw[-] (11.5,-3.25) -- (otimes);

    \node[rotate=0,minimum width=0.5cm] (u) at (13.5,-6.75) {\small $\widetilde{\mathbf{u}}$};

    \node[rotate=0,minimum width=0.5cm] (loss) at (13.5,-9) {\small $\mathcal{L} \left( \mathbf{u}_\text{pred}, \widetilde{\mathbf{u}} \right)$};
	
	\draw[->] (otimes) -- (u);
	\draw[->] (u) -- (loss);

	\end{tikzpicture}
	\caption{PTFONet architecture; the shown text blocks below the layers describe their output tensor dimensions. The symbol $\odot$ and $\otimes$ are, respectively, the Hadamard and tensor products, as explained in Section \ref{sec:ptfonet}. $\mathcal{L}$ is the relative $\mathcal{L}^2-$error loss function between the reference (true) displacement $\mathbf{u}_\text{pred}$ and the predicted displacement $\widetilde{\mathbf{u}}$, given by Eq. \ref{eq:loss}. ReLU activation functions are used after each block with the exception of the last MLP blocks in the branch and trunk networks, where no activation function was used; a tanh activation function was used at the end of the last MLP block in the BC network.}
	\label{fig:ptfonet_architecture}
\end{figure}
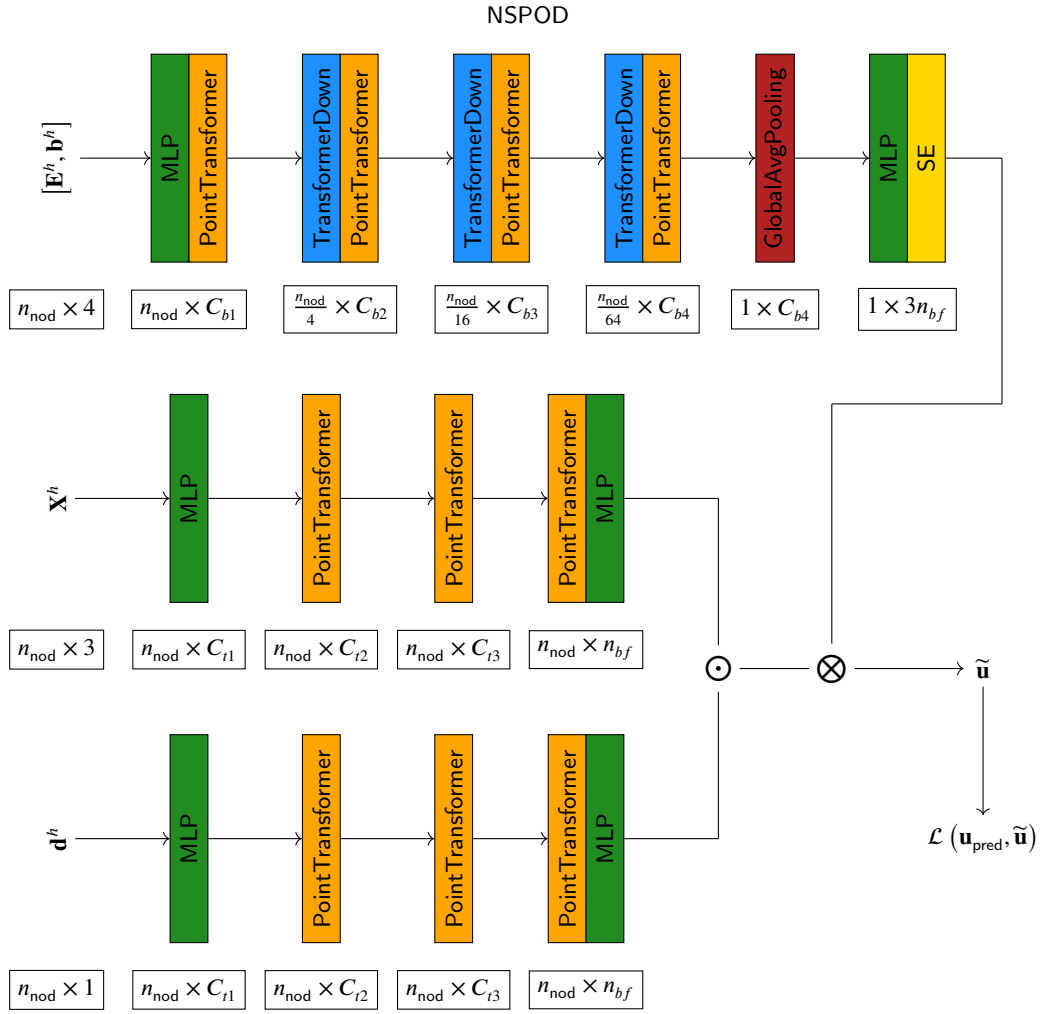

\noindent The overall architecture of PTFONet is shown in Fig. \ref{fig:ptfonet_architecture}. Edge indices $\mathcal{E}_i$ and nodal coordinates $\mathbf{X}^h_i$ are also provided to the PointTransformer layers in order to compute nodal neighbors. Normally the edge indices would not need to be provided and the appropriate connectivity would be calculated via $k-$nearest neighbors ($k$NN) \cite{PointTransformer} but we provide the edge indices to all the subnetworks of PTFONet so that we do not waste computational resources on unnecessary computation of neighbors; however, the edge indices are still recomputed via $k$NN after each pooling operation in the branch network. Here it is also important to point out that in order for our precomputed edge indices to perform similarly to those computed by $k$NN we require the mesh to be of good quality, meaning that for all the nodes $\mathbf{X}^h_i \in \mathcal{T}^h$, there is some $\varepsilon > 0$ such that the ball of radius $\varepsilon$ centered around $\mathbf{X}^h_i$ contains all nodes $\mathbf{X}^h_j \in \mathcal{N}_{\mathcal{E}_i} \left( \mathbf{X}^h_i \right)$ and contains no nodes $\mathbf{X}^h_k \notin \mathcal{N}_{\mathcal{E}_i} \left( \mathbf{X}^h_i \right)$, where $\mathcal{N}_{\mathcal{E}_i} \left( \mathbf{X}^h_i \right)$ is the neighborhood of $\mathbf{X}^h_i$ as defined by the edge indices $\mathcal{E}_i$.

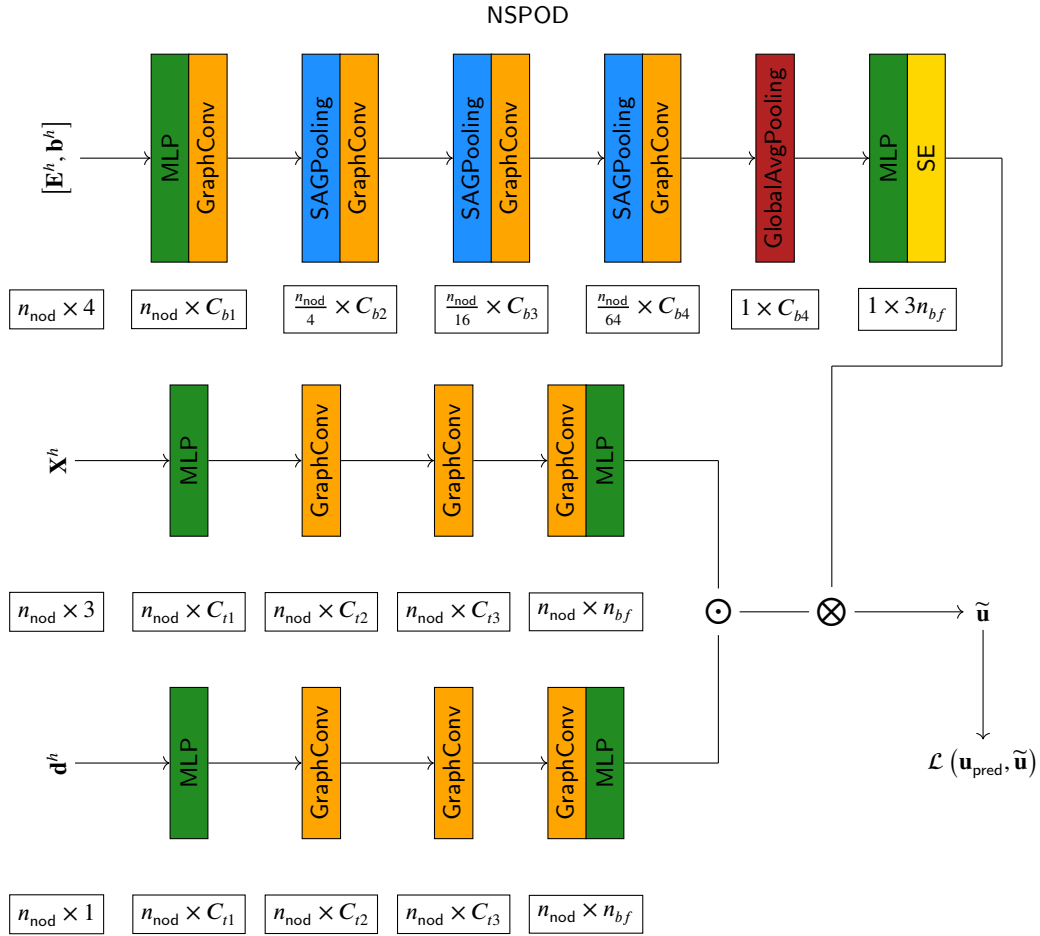
\begin{figure}
\centering
\begin{tikzpicture}

	\node[rotate=90,minimum width=2.75cm] (x1) at (1.25,0) {\small$\left[ \mathbf{E}^h, \mathbf{b}^h \right]$};
	\node[h,rotate=90,minimum width=2.75cm] (mlp1) at (2.75,0) {\small MLP};
	\node[conv,rotate=90,minimum width=2.75cm] (pntb1) at (3.25,0) {\small GraphConv};
	\node[fc,rotate=90,minimum width=2.75cm] (td1) at (4.75,0) {\small SAGPooling};
	\node[conv,rotate=90,minimum width=2.75cm] (pntb2) at (5.25,0) {\small GraphConv};
	\node[fc,rotate=90,minimum width=2.75cm] (td2) at (6.75,0) {\small SAGPooling};
	\node[conv,rotate=90,minimum width=2.75cm] (pntb3) at (7.25,0) {\small GraphConv};
	\node[fc,rotate=90,minimum width=2.75cm] (td3) at (8.75,0) {\small SAGPooling};
	\node[conv,rotate=90,minimum width=2.75cm] (pntb4) at (9.25,0) {\small GraphConv};
	\node[pool,rotate=90,minimum width=2.75cm] (pool) at (10.75,0) {\small GlobalAvgPooling};
    \node[h,rotate=90,minimum width=2.75cm] (mlp2) at (12.25,0) {\small MLP};
    \node[bn,rotate=90,minimum width=2.75cm] (SE) at (12.75,0) {\small SE};

    \node[view,rotate=0,minimum width=1cm] (b1) at (1.25,-2) {\small $n_\text{nod} \times 4$};
    \node[view,rotate=0,minimum width=1cm] (b2) at (3,-2) {\small $n_\text{nod} \times C_{b1}$};
    \node[view,rotate=0,minimum width=1cm] (b4) at (5,-2) {\small $\frac{n_\text{nod}}{4} \times C_{b2}$};
    \node[view,rotate=0,minimum width=1cm] (b5) at (7,-2) {\small $\frac{n_\text{nod}}{16} \times C_{b3}$};
    \node[view,rotate=0,minimum width=1cm] (b6) at (9,-2) {\small $\frac{n_\text{nod}}{64} \times C_{b4}$};
    \node[view,rotate=0,minimum width=1cm] (b7) at (10.75,-2) {\small $1 \times C_{b4}$};
    \node[view,rotate=0,minimum width=1cm] (b8) at (12.5,-2) {\small $1 \times 3 n_{bf}$};

    \node[rotate=90,minimum width=2cm] (x2) at (1.25,-4) {\small$\mathbf{X}^h$};
    \node[h,rotate=90,minimum width=2cm] (mlp3) at (3,-4) {\small MLP};
	\node[conv,rotate=90,minimum width=2cm] (pntb5) at (4.75,-4) {\small GraphConv};
    \node[conv,rotate=90,minimum width=2cm] (pntb6) at (6.5,-4) {\small GraphConv};
    \node[conv,rotate=90,minimum width=2cm] (pntb7) at (8,-4) {\small GraphConv};
    \node[h,rotate=90,minimum width=2cm] (mlp4) at (8.5,-4) {\small MLP};

    \node[view,rotate=0,minimum width=1cm] (t1) at (1.25,-6) {\small $n_\text{nod} \times 3$};
    \node[view,rotate=0,minimum width=1cm] (t2) at (3,-6) {\small $n_\text{nod} \times C_{t1}$};
    \node[view,rotate=0,minimum width=1cm] (t3) at (4.75,-6) {\small $n_\text{nod} \times C_{t2}$};
    \node[view,rotate=0,minimum width=1cm] (t4) at (6.5,-6) {\small $n_\text{nod} \times C_{t3}$};
    \node[view,rotate=0,minimum width=1cm] (t5) at (8.25,-6) {\small $n_\text{nod} \times n_{bf}$};

    \node[rotate=90,minimum width=2cm] (x3) at (1.25,-8) {\small$\mathbf{d}^h$};
    \node[h,rotate=90,minimum width=2cm] (mlp5) at (3,-8) {\small MLP};
	\node[conv,rotate=90,minimum width=2cm] (pntb8) at (4.75,-8) {\small GraphConv};
    \node[conv,rotate=90,minimum width=2cm] (pntb9) at (6.5,-8) {\small GraphConv};
    \node[conv,rotate=90,minimum width=2cm] (pntb10) at (8,-8) {\small GraphConv};
    \node[h,rotate=90,minimum width=2cm] (mlp6) at (8.5,-8) {\small MLP};

    \node[view,rotate=0,minimum width=1cm] (t1) at (1.25,-10) {\small $n_\text{nod} \times 1$};
    \node[view,rotate=0,minimum width=1cm] (t2) at (3,-10) {\small $n_\text{nod} \times C_{t1}$};
    \node[view,rotate=0,minimum width=1cm] (t3) at (4.75,-10) {\small $n_\text{nod} \times C_{t2}$};
    \node[view,rotate=0,minimum width=1cm] (t4) at (6.5,-10) {\small $n_\text{nod} \times C_{t3}$};
    \node[view,rotate=0,minimum width=1cm] (t5) at (8.25,-10) {\small $n_\text{nod} \times n_{bf}$};
		
	\draw[->] (x1) -- (mlp1);
	\draw[->] (pntb1) -- (td1);
	\draw[->] (pntb2) -- (td2);
	\draw[->] (pntb3) -- (td3);
    \draw[->] (pntb4) -- (pool);
    \draw[->] (pool) -- (mlp2);
	
	\draw[->] (x2) -- (mlp3);
	\draw[->] (mlp3) -- (pntb5);
	\draw[->] (pntb5) -- (pntb6);
	\draw[->] (pntb6) -- (pntb7);
	
	\draw[->] (x3) -- (mlp5);
	\draw[->] (mlp5) -- (pntb8);
	\draw[->] (pntb8) -- (pntb9);
	\draw[->] (pntb9) -- (pntb10);

    \node[rotate=90,minimum width=0.5cm] (odot) at (10,-6) {\Large $\odot$};
    \node[rotate=90,minimum width=0.5cm] (otimes) at (11.5,-6) {\Large $\otimes$};
        
	\draw[-] (SE) -- (13.75,0);
	\draw[-] (13.75,0) -- (13.75,-2.75);
	
	\draw[-] (mlp4) -- (10,-4);
    \draw[-] (mlp6) -- (10,-8);
	\draw[-] (10,-4) -- (odot);
    \draw[-] (odot) -- (10,-8);
    \draw[-] (odot) -- (otimes);
	\draw[-] (13.75,-2.75) -- (11.5,-2.75);
	\draw[-] (11.5,-2.75) -- (otimes);

    \node[rotate=0,minimum width=0.5cm] (u) at (13.5,-6) {\small $\widetilde{\mathbf{u}}$};

    \node[rotate=0,minimum width=0.5cm] (loss) at (13.5,-8) {\small $\mathcal{L} \left( \mathbf{u}_\text{pred}, \widetilde{\mathbf{u}} \right)$};
	
	\draw[->] (otimes) -- (u);
	\draw[->] (u) -- (loss);

	\end{tikzpicture}
	\caption{Graph-based architecture; the shown text blocks below the layers describe their output tensor dimensions. The symbol $\odot$ and $\otimes$ are, respectively, the Hadamard and tensor products, as explained in Section \ref{sec:ptfonet}. $\mathcal{L}$ is the relative $\mathcal{L}^2-$error loss function between the reference (true) displacement $\mathbf{u}_\text{pred}$ and the predicted displacement $\widetilde{\mathbf{u}}$, given by Eq. \ref{eq:loss}. ReLU activation functions are used after each block with the exception of the last MLP blocks in the branch and trunk networks, where no activation function was used; a tanh activation function was used at the end of the last MLP block in the BC network.}
	\label{fig:gcn_architecture}
\end{figure}

We include an ablation study in Section \ref{sec:ablation_study}, where we substitute: PointTransformer convolution layers for GraphConv \cite{graphconv} convolution layers, and Farthest Point Sampling (FPS) \cite{fps, PointTransformer} and kNN-based pooling layers for SAGPooling-based \cite{sagpooling1, sagpooling2} pooling layers. This creates an analogous graph-based DeepONet, which we call GraphONet, and allows for a straightforward comparison of performance between PointTransformer and a \textit{state-of-the-art} graph convolution layer, GraphConv. We choose these specific graph-based convolution and pooling layers because they achieved the maximum accuracy in our numerical experiments with other configurations of graph-based DeepONets (results not shown). The architecture of this graph-based DeepONet is shown in Fig. \ref{fig:gcn_architecture}.

%%%%%%%%%%%%%%%%%%%%%%%%%%%%%%%%%%%%%%%%%%%

\subsection{Training methodology} \label{sec:training}
\noindent To train PTFONet we build a dataset $D$ consisting of $N_s$ samples given by
\begin{equation} \label{eq:dataset}
D = \left\{ \left( \mathcal{T}^h_i, \mathbf{E}^h_i, \mathbf{b}^h_i, \mathbf{X}^h_i, \mathbf{d}^h_i, \mathcal{E}_i, \mathbf{u}_{\text{pred}, i} \right) \right\}_{i = 1}^{N_s}.
\end{equation}

\noindent Naturally, the discretized domain $\mathcal{T}^h_i$ determines the matrix of coordinates $\mathbf{X}^h_i$ and set of edge indices $\mathcal{E}_i$; $\mathbf{E}^h_i$, $\mathbf{b}^h_i$, $\mathbf{d}^h_i$ and $\mathbf{u}^h_{\text{pred}, i}$ here respectively denote the vectors/matrices of nodal values of Youngs' modulus, body forces, locations where homogeneous Dirichlet BCs have been applied, and reference displacement solution, all computed on $\mathcal{T}^h_i$. The training process of PTFONet is performed by minimizing the nodal-wise relative error between the inferred and reference displacement solutions, and is given by
\begin{equation} \label{eq:loss}
\text{arg\,min}_{\pmb{\theta}} \frac{1}{N_s} \sum_{i = 1}^{N_s} \sqrt{\frac{\left( \mathcal{G}_\text{PTF} \left( \mathbf{E}^h_i, \mathbf{b}^h_i, \mathbf{X}^h_i, \mathcal{E}_i \right) \left( \mathbf{X}^h_i, \mathbf{d}^h_i, \mathcal{E}_i \right) - \mathbf{u}_{\text{pred}, i} \right)^\text{T} \left( \mathcal{G}_\text{PTF} \left( \mathbf{E}^h_i, \mathbf{b}^h_i, \mathbf{X}^h_i, \mathcal{E}_i \right) \left( \mathbf{X}^h_i, \mathbf{d}^h_i, \mathcal{E}_i \right) - \mathbf{u}_{\text{pred}, i} \right)}{\mathbf{u}_\text{pred}^\text{T} \mathbf{u}_\text{pred}}},
\end{equation}

\noindent where $\pmb{\theta}$ denote all the parameters in PTFONet. We also use the same relative error function as our accuracy function.

All our discretized meshes $\mathcal{T}^h_i$ started as CAD geometries; we use a total of four meshes which were obtained from the ABC dataset \cite{ABCCAD}. We then mesh these geometries with $P^1$ conforming finite elements with Gmsh \cite{gmsh} and process them with FEniCSx \cite{DOLFINx} in order to assign the appropriate BCs, as we are describing in the next paragraph. We ensure that the chosen geometries have around $10,000$ nodes in order to not require excessive computational power during training. We also use FEniCSx and the sparse direct solver MUMPS \cite{MUMPS} to solve the corresponding linear elasticity systems on these geometries in order to obtain $\mathbf{u}_\text{pred}$.

As previously mentioned, in this study we use spatially constant $E$ and $\mathbf{b}$ for each considered sample in the dataset $D$. We obtain $E$ for a specific sample by generating values $E \sim N \left( 1, 1 \right)$ until $E > E_\text{min}$, where $E_\text{min} = 0.3$. To obtain $\mathbf{b}$ for a specific sample, we generate it as in $\mathbf{b} \sim N \left( \mathbf{0}_3, \mathbf{I}_3 \right)$, where $\mathbf{0}_3$ and $\mathbf{I}_3$ are the three-dimensional zero vector and identity matrix, respectively. $N \left( \pmb{\mu}, \pmb{\Sigma} \right)$ denotes a normal distribution with mean vector $\pmb{\mu}$ and covariance matrix $\pmb{\Sigma}$. Note that despite our use of spatially constant $E$, we still construct a matrix $\mathbf{E}^h$ with repeated numbers for each node. Furthermore, a matrix $\mathbf{b}^{h}$ is constructed by applying $L^{2}$ projection to the right-hand side vector $\mathbf{f}$ in Eq \eqref{eq:linear_system}, which leads a non-constant matrix. This is to allow for a straightforward extension to when spatially varying functions are used as inputs. To obtain $d$ we first sample $p$ from a uniform distribution, namely $p \sim U \left( 0.1, 0.5 \right)$. This represents the maximum proportion of boundary that gets homogeneous Dirichlet BC assigned, i.e. $0 < \frac{\left| \partial\Omega_0^D \right|}{\left| \partial\Omega_0 \right|} \leq p$. Then create two lists: list $A$, which contains all the boundary surfaces for each of the considered geometries, randomly permuted; and list $B$, which will contain the boundary surfaces that form $\partial\Omega_0^D$. We go through the list of surfaces in list $A$ exhaustively, adding these to list $B$ whilst $\frac{\left| \partial\Omega_0^D \right|}{\left| \partial\Omega_0 \right|} \leq p$ holds, and skipping any surfaces that would lead to $\frac{\left| \partial\Omega_0^D \right|}{\left| \partial\Omega_0 \right|} > p$.

We implemented all the networks in this study using PyTorch \cite{pytorch} and trained using the AdamW optimizer \cite{AdamW}, with a batch size of 10, a learning rate of $10^{-4}$, and weight decay of $5 \times 10^{-5}$. The training process stops if the average relative error (Eq. \ref{eq:loss}) of the test samples does not decrease for 100 consecutive epochs. In all our numerical experiments, the Krylov-based iterative solvers terminate when either of the two following criteria are satisfied: $\| \mathbf{r}^{\left( i \right)} \|_2 \leq 10^{-12},$ or $\frac{\| \mathbf{r}^{\left( i \right)} \|_2}{\| \mathbf{r}^{\left( 0 \right)} \|_2} \leq 10^{-8}$. We performed all our network training and data generation on a computational cluster with an Intel Xeon Gold 6542Y and an NVIDIA H100 NVL Tensor Core GPU. Other experiments were conducted using computational resources and services at the Center for Computation and Visualization, Brown University, where each computing node is equipped with an AMD EPYC 9554 64-Core Processor (256GB) and an NVIDIA L40S GPU (48GB).

%%%%%%%%%%%%%%%%%%%%%%%%%%%%%%%%%%%%%%%%%%%

\subsection{Neural Proper Orthogonal Decomposition} \label{sec:neural_POD}
\begin{figure}
\centering
    \begin{tikzpicture}[very thick, fill opacity=1, text opacity=1]

    \draw[fill=blue!50, fill opacity=0.5, rounded corners] (0, 0) rectangle (2.5, 0.8) node[midway]{$\mathcal{G}_\text{PTF} \left( \mathbf{f}_1 \right) \left( \mathbf{X}_1 \right)$};
    \draw[fill=blue!50, fill opacity=0.5, rounded corners] (0, -1) rectangle (2.5, -0.2) node[midway]{$\mathcal{G}_\text{PTF} \left( \mathbf{f}_2 \right) \left( \mathbf{X}_2 \right)$};
    \draw[fill=blue!50, fill opacity=0.5, rounded corners] (0, -2) rectangle (2.5, -1.2) node[midway]{$\mathcal{G}_\text{PTF} \left( \mathbf{f}_3 \right) \left( \mathbf{X}_3 \right)$};
    \draw[fill=blue!50, fill opacity=0.5, rounded corners] (0, -3) rectangle (2.5, -2.2) node[midway]{$\dots$};
    \draw[fill=blue!50, fill opacity=0.5, rounded corners] (0, -4) rectangle (2.5, -3.2) node[midway]{$\mathcal{G}_\text{PTF} \left( \mathbf{f}_m \right) \left( \mathbf{X}_m \right)$};

    \draw [decorate, decoration={brace, amplitude=5pt}] (2.6, 0.8) -- (2.6, -4);

    \draw[fill=green!50, fill opacity=0.5, rounded corners] (3.5, -2.7) rectangle (10.25, -0.5) node[midway]{$\mathbf{U} = \begin{bmatrix}
    \mathcal{G}_\text{PTF} \left( \mathbf{f}_1 \right) \left( \mathbf{X}_{1,1} \right) & \dots & \mathcal{G}_\text{PTF} \left( \mathbf{f}_1 \right) \left( \mathbf{X}_{1,n_\text{nod}} \right) \\
    \mathcal{G}_\text{PTF} \left( \mathbf{f}_2 \right) \left( \mathbf{X}_{2,1} \right) & \dots & \mathcal{G}_\text{PTF} \left( \mathbf{f}_2 \right) \left( \mathbf{X}_{2,n_\text{nod}} \right) \\
    \vdots & \ddots & \vdots \\
    \mathcal{G}_\text{PTF} \left( \mathbf{f}_m \right) \left( \mathbf{X}_{m,1} \right) & \dots & \mathcal{G}_\text{PTF} \left( \mathbf{f}_m \right) \left( \mathbf{X}_{m,n_\text{nod}} \right)
    \end{bmatrix}$};

    \draw [arrows = {-Latex[width'=0pt .5, length=10pt]}] (2.75,-1.6) -- (3.5,-1.6);

    \draw[fill=red!50, fill opacity=0.5, rounded corners] (11, -2.7) rectangle (16.25, -0.5) node[midway]{$\begin{aligned}
        &\pmb{\Sigma} = \frac{1}{m - 1}\mathbf{U} \mathbf{U}^\text{T} \\
        &\text{POD: } \left\{ \left\{ \lambda_1, \lambda_2, \dots \right\}, \left\{ \pmb{\phi}_1, \pmb{\phi}_2, \dots \right\} \right\} \\
        &\mathbf{P}^\text{T} = \begin{bmatrix}
            \pmb{\phi}_1 & \pmb{\phi}_2 & \dots & \pmb{\phi}_m
        \end{bmatrix}
    \end{aligned}$};

    \draw [arrows = {-Latex[width'=0pt .5, length=10pt]}] (10.25,-1.6) -- (11,-1.6);
    
    \end{tikzpicture}
	\caption{Pipeline for NSPOD creation. From left to right: inferring $m$ solutions for $m$ different input combinations $\left\{ \mathbf{f}_{i} \right\}_{i=1}^{m}$, defined in~\eqref{eq:noised_input}, assembly of the $\mathbf{U}$ matrix, creation of the appropriate covariance matrix $\pmb{\Sigma}$ and its spectral decomposition, final assembly of the prolongation matrix $\mathbf{P}$. $\mathbf{X}_{i,j}$ here denote the $j-$th node for the $i-$th modified right-hand side vector.}
	\label{fig:NeuralPOD}
\end{figure}
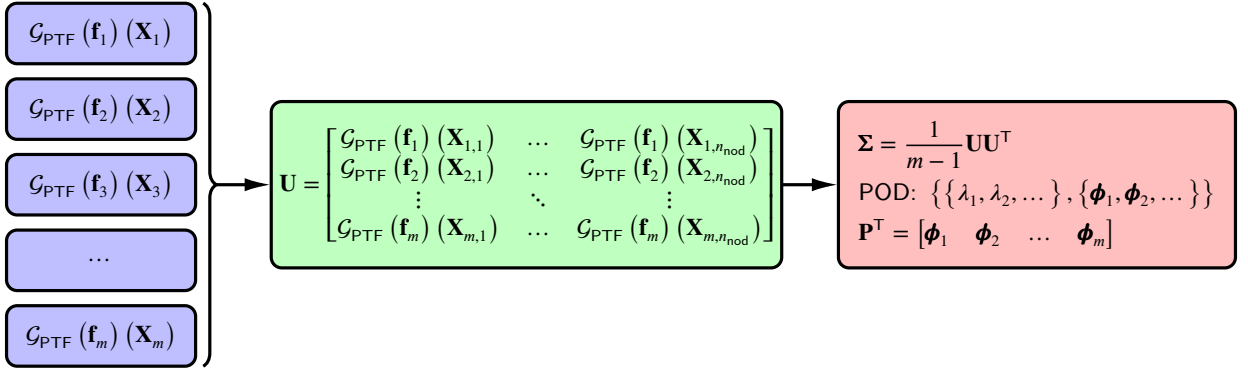

\noindent In this study we introduce the NSPOD preconditioner for Krylov-based methods, which is based on creating a prolongation operator $\mathbf{P}$ in Eq. \ref{eq:krylov_hybrid_preconditioner} from a POD of the solution space spanned by applying different input functions. To be more precise, we first train PTFONet with the desired combination of geometries, each sample having different $E$, $\mathbf{b}$, and $d$. Then we use the trained PTFONet to infer $\widetilde{\mathbf{u}}$ on a number $m$ of unseen cases, which are generated not by varying physical parameters, but by introducing random noise to the right-hand side vector $\mathbf{f}$ to form a neighborhood of the original problem. In particular, the $i$-th modified right-hand side vector is defined as 
\begin{equation}\label{eq:noised_input}
    \mathbf{f}_{i} = \mathbf{f}+\eta \mathbf{s}_{i},
\end{equation}
where $\eta$ is a scaling factor and each component of the noise vector $\mathbf{s}_{i}$ is sampled from a uniform distribution $\mathcal{U}\left(\mathbf{f}_{\text{min}}, \mathbf{f}_{\text{max}}\right)$.
Here, $\mathbf{f}_{\text{min}}$ and $\mathbf{f}_{\text{max}}$ represent the minimum and maximum elements of the nominal vector $\mathbf{f}$, respectively.
Note that we set $\eta=0.01$ in this paper. We then build a matrix $\mathbf{U}$ where the rows correspond to these inferred $\widetilde{\mathbf{u}}$. Then we build the corresponding covariance matrix $\pmb{\Sigma} = \frac{1}{m-1} \mathbf{U} \mathbf{U}^\text{T}$ and perform its spectral decomposition with eigenvectors $\pmb{\phi}_1, \pmb{\phi}, \dots, \pmb{\phi}_m$. We finally build the desired prolongation operator as $\mathbf{P}^\text{T} = \begin{bmatrix} \pmb{\phi}_1 & \pmb{\phi}_2 & \dots & \pmb{\phi}_m \end{bmatrix}$. The process is depicted in Fig. \ref{fig:NeuralPOD}.

%%%%%%%%%%%%%%%%%%%%%%%%%%%%%%%%%%%%%%%%%%%

\section{Results} \label{sec:results}
\noindent In this section we describe a series of numerical experiments and report their results. In Section \ref{sec:single_geometry} we show the performance of the NSPOD preconditioner for Krylov-based iterative linear solvers when the associated PTFONet is trained on one geometry. In Section \ref{sec:several_geometries} we evaluate the performance of the NSPOD preconditioner for Krylov-based iterative linear solvers when the associated PTFONet is trained on several geometries. In Section \ref{sec:ablation_study} we show an ablation study which compared the performance of PTFONet against the performance of an analogous graph-based DeepONet. We finalize the Results section by reporting the computational cost of the NSPOD preconditioner in Section \ref{sec:computational_cost}. Since we are using PTFONets trained on different datasets we summarize this information, including the names, number of train/test samples, and brief description of the network/dataset, in Table \ref{tab:datasets}.
\begin{table}
	\centering
	\caption{Names of the different versions of PTFONet used, number of samples, and brief descriptions of the considered datasets.} \label{tab:datasets}
	\begin{tabular}{ccc}
		\toprule
		\textbf{Name} & \textbf{Number of samples (train/test)} & \textbf{Brief description of network/dataset} \\
		\midrule
		PTFONetV1 & 4000 (3600/400) & One geometry \\
        PTFONetV2 & 6000 (5400/600) & Three geometries $-$ no unseen geometries in the test set \\
        PTFONetV3 & 4000 (3600/400) & One geometry $-$ used to infer on unseen geometries \\
        GraphONet & 6000 (5400/600) & Graph-based DeepONet $-$ same dataset as PTFONetV2 \\
		\bottomrule
	\end{tabular}
\end{table}

%%%%%%%%%%%%%%%%%%%%%%%%%%%%%%%%%%%%%%%%%%%

\subsection{Training on a single geometry} \label{sec:single_geometry}
\begin{figure}
	\centering
	\subfloat[Geometry 1]{
		\includegraphics[width=0.23\linewidth]{Figure_5a.png}
	}
	\subfloat[Geometry 2]{
		\includegraphics[width=0.329\linewidth]{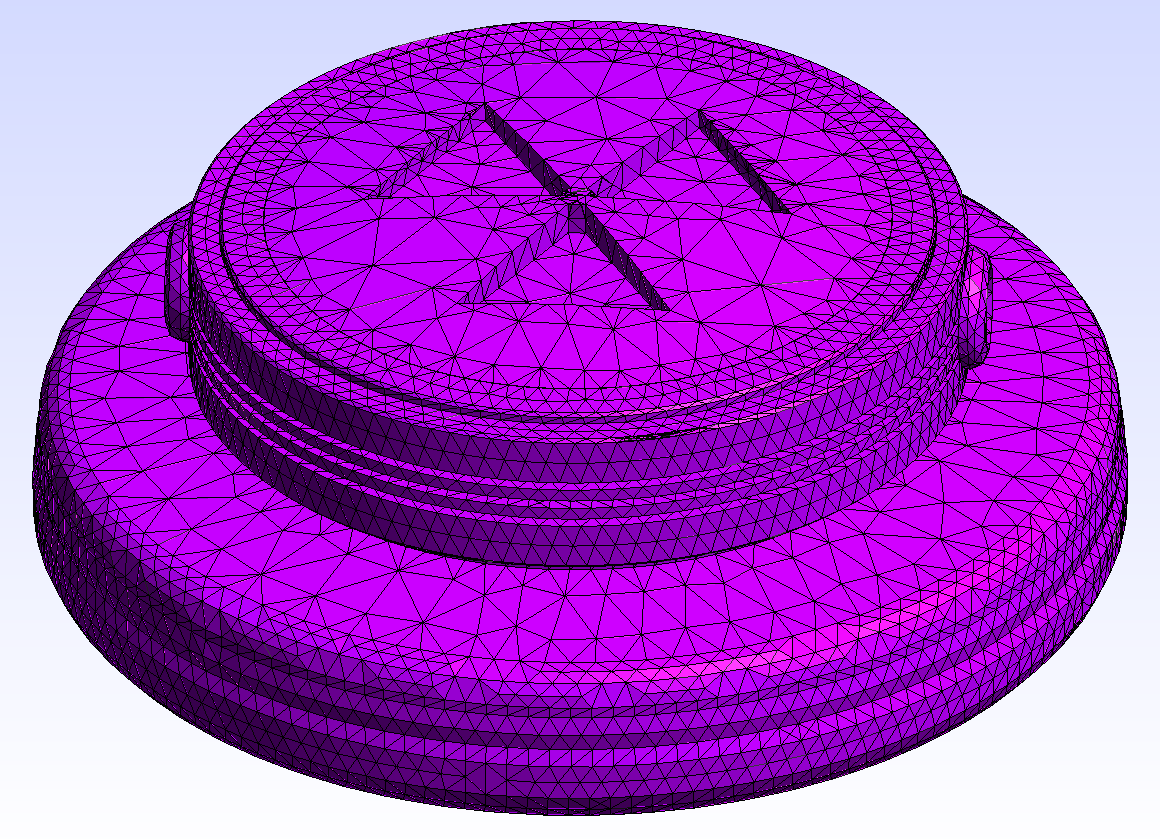}
	}
	\subfloat[Geometry 3]{
		\includegraphics[width=0.4\linewidth]{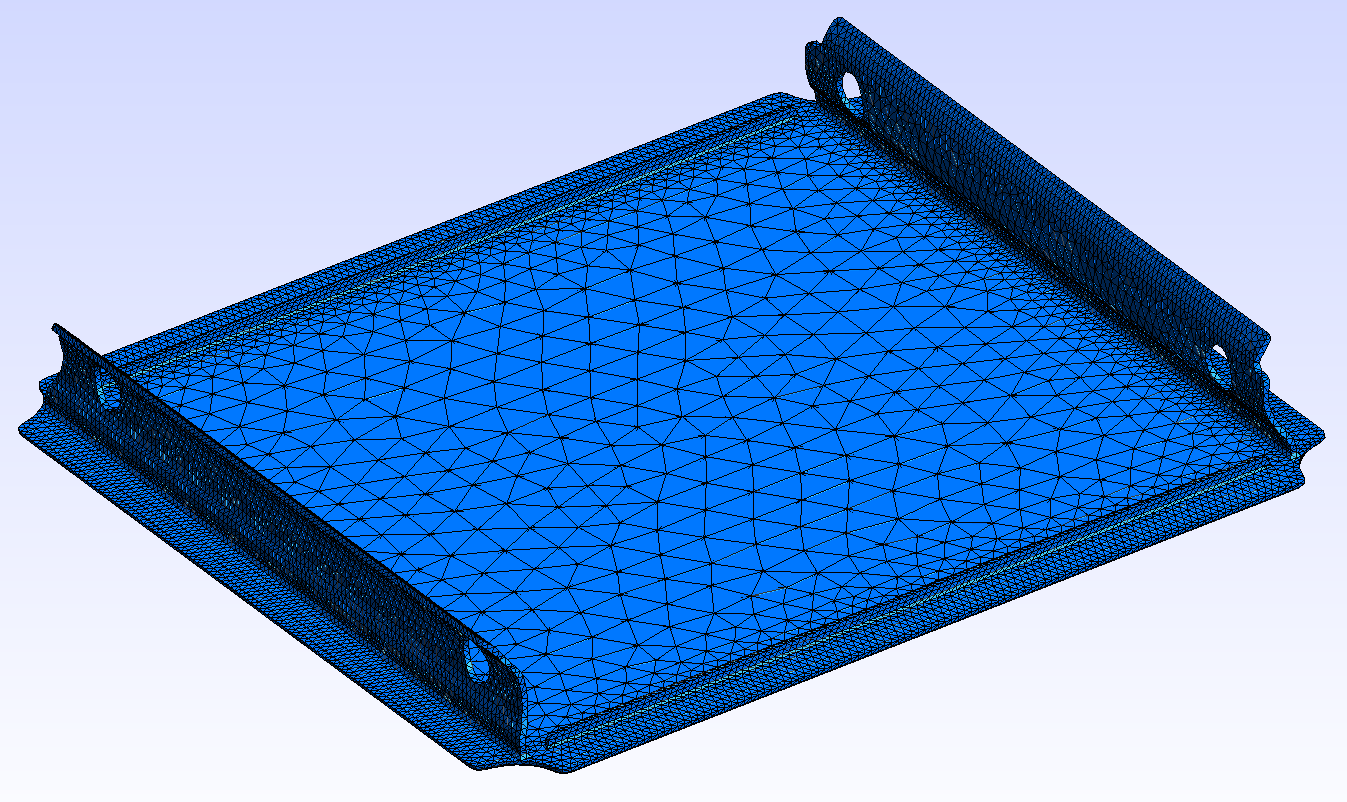}
	}

	\caption{Discretized geometries for PTFONetV1 training on a single geometry. $P^{1}$ conforming tetrahedra finite elements are used for the discretization.
		\textbf{(a)} The mesh of Geometry 1 consists of $11,411$ nodes and $48,215$ elements.
		\textbf{(b)} The mesh of Geometry 2 consists of $11,163$ nodes and $43,216$ elements. 
		\textbf{(c)} The mesh of Geometry 3 consists of $15,563$ nodes and $45,258$ elements. 
	}\label{fig:geometries}
\end{figure}

\noindent We trained PTFONet three times with a different single geometry each time; these geometries are shown in Fig. \ref{fig:geometries}. The employed architecture, PTFONetV1 from now onwards, is shown in Fig. \ref{fig:ptfonet_architecture}; the parameters shown in the figure are described in Table \ref{tab:channels}. We employed three datasets each containing $N_s = 4000$ samples of the same geometry but different combinations of $E$, $\mathbf{b}$, and $d$. The relative errors (Eq. \ref{eq:loss}) on the respective test sets achieved by the trained PTFONetV1 are $12.14\%$ for Geometry 1, $13.50\%$ for Geometry 2, and $17.15\%$ for Geometry 3.

\begin{table}
	\centering
	\caption{Number of channels for the layers shown in Fig. \ref{fig:ptfonet_architecture} for PTFONetV1.} \label{tab:channels}
	\begin{tabular}{cc}
		\toprule
		\textbf{Parameter} & \textbf{Number of channels} \\
		\midrule
		$C_\text{b1}$ & 48 \\
        $C_\text{b2}$ & 96 \\
        $C_\text{b3}$ & 192 \\
        $C_\text{b4}$ & 384 \\
        $C_\text{t1}$ & 48 \\
        $C_\text{t2}$ & 96 \\
        $C_\text{t3}$ & 192 \\
        $N_\text{bf}$ & 256 \\
		\bottomrule
	\end{tabular}
\end{table}

\begin{figure}
	\centering
	\includegraphics[width=1\linewidth, clip]{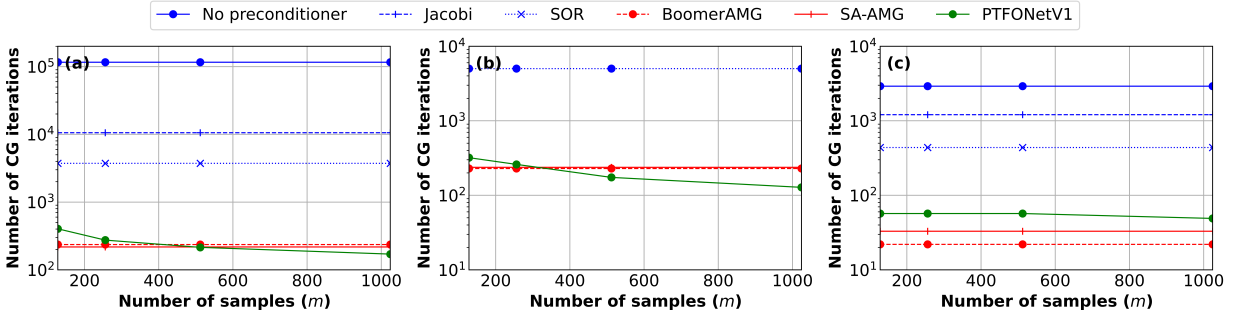}
	\caption{Number of CG iterations \textit{vs} number of samples employed to create the NSPOD preconditioner for the three considered geometries; we employed PTFONetV1 in this scenario. We consider $m = 128, 256, 512, \text{ and } 1024$.
    \textbf{(a)} Geometry 1. 
    \textbf{(b)} Geometry 2, where CG with no preconditioner and/or with Jacobi preconditioner did not converge, thus the corresponding results are not shown.
    \textbf{(c)} Geometry 3.} \label{fig:num_iterations_single_geometry}
\end{figure}

We now compare the number of iterations needed for convergence of different preconditioners when a CG iterative solver is used. Figure \ref{fig:num_iterations_single_geometry} shows the number of iterations to convergence of no preconditioner, a Jacobi preconditioner, SOR preconditioner, BoomerAMG \cite{hypre}, SA-AMG~\cite{van2001convergence} from PETSc preconditioner, and our NSPOD preconditioner for different number of samples $m$. BoomerAMG is a classical Ruge-St\"uben style AMG. Specifically, it uses Falgout coarsening, canonical interpolation, and Symmetric SOR/Jacobi smoother. On the other hand, PETSc's SA-AMG uses aggregation coarsening, smoothed interpolation with damped Jacobi, and Chebyshev smoother. These preconditioners are tested on an unseen combination of $E$, $\mathbf{b}$, and $d$ for geometries 1, 2, and 3. Figure \ref{fig:num_iterations_single_geometry} shows that the performance of the NSPOD preconditioner using PTFONetV1 outperforms even \text{state-of-the-art} preconditioners such as BoomerAMG or PETSc's SA-AMG when the number of samples $m$ to construct the POD is large enough. For Geometry 1, we have that for $m = 1,024$ the number of iterations for NSPOD with PTFONetV1 is 171 whilst for BoomerAMG is 218 and for PETSc's SA-AMG is 237. For Geometry 2, we have that for $m = 1,024$ the number of iterations for NSPOD with PTFONetV1 is 128, respectively, whilst for BoomerAMG is 238 and for PETSc's SA-AMG is 229. For Geometry 3 we have a different trend, where the performance of NSPOD is poorer than that of the standard AMG preconditioners evaluated; we have that for $m = 1,024$ the number of iterations for NSPOD with PTFONetV1 is 49 whilst for BoomerAMG is 33 and for PETSc's SA-AMG is 22.

%%%%%%%%%%%%%%%%%%%%%%%%%%%%%%%%%%%%%%%%%%%

\subsection{Training on several geometries} \label{sec:several_geometries}
\noindent We now trained PTFONet on the geometries shown in Fig. \ref{fig:geometries} at the same time in order to check the potential geometric transferability of PTFONet. The employed architecture, PTFONetV2 from now onwards, is the same as in the previous section, shown in Fig. \ref{fig:ptfonet_architecture}; however the parameters are different, now described in Table \ref{tab:channels_several}. Since we wanted to encode several geometries in the same network we increased the capacity of the network. Also, we increased the value of weight decay to $10^{-2}$ in order to prevent potential overfitting caused by this substantial increase in capacity. We employed a dataset containing $N_s = 6000$ samples, $2000$ of each of the considered geometries, and with different combinations of $E$, $\mathbf{b}$, and $d$; both the train and test sets of this dataset contain the three geometries so there is no unseen geometry in the test set. The relative error on the test set achieved by the trained PTFONetV2 is $16.91\%$.
\begin{table}
	\centering
	\caption{Number of channels for the layers shown in Fig. \ref{fig:ptfonet_architecture} for PTFONetV2.} \label{tab:channels_several}
	\begin{tabular}{cc}
		\toprule
		\textbf{Parameter} & \textbf{Number of channels} \\
		\midrule
		$C_\text{b1}$ & 64 \\
        $C_\text{b2}$ & 128 \\
        $C_\text{b3}$ & 256 \\
        $C_\text{b4}$ & 768 \\
        $C_\text{t1}$ & 64 \\
        $C_\text{t2}$ & 128 \\
        $C_\text{t3}$ & 256 \\
        $N_\text{bf}$ & 512 \\
		\bottomrule
	\end{tabular}
\end{table}

\begin{figure}
	\centering
	\includegraphics[width=1\linewidth, clip]{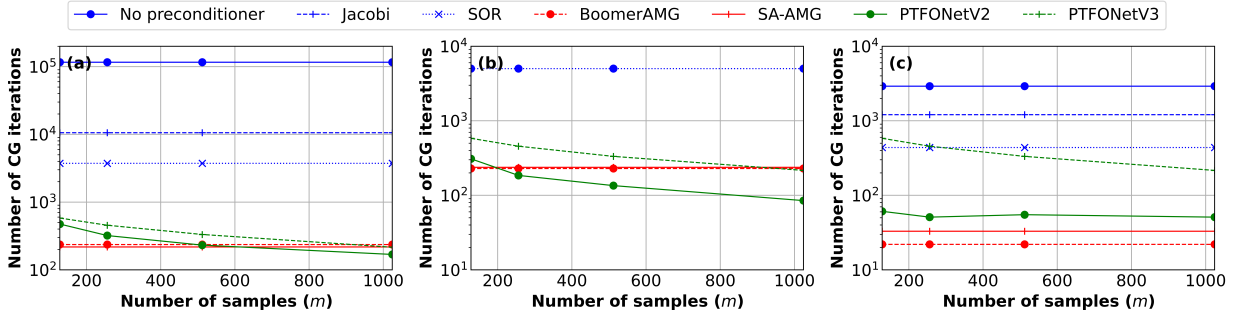}
	\caption{Number of CG iterations \textit{vs} number of samples employed to create the NSPOD preconditioner for the three considered geometries; we employed PTFONetV2 and PTFONetV3 in this scenario. We consider $m = 128, 256, 512, \text{ and } 1024$.
    \textbf{(a)} Geometry 1. 
    \textbf{(b)} Geometry 2, where CG with no preconditioner and/or with Jacobi preconditioner did not converge, thus the corresponding results are not shown.
    \textbf{(c)} Geometry 3.} \label{fig:num_iterations_several_geometries}
\end{figure}

\begin{figure}
	\centering
	\includegraphics[width=0.4\linewidth, clip]{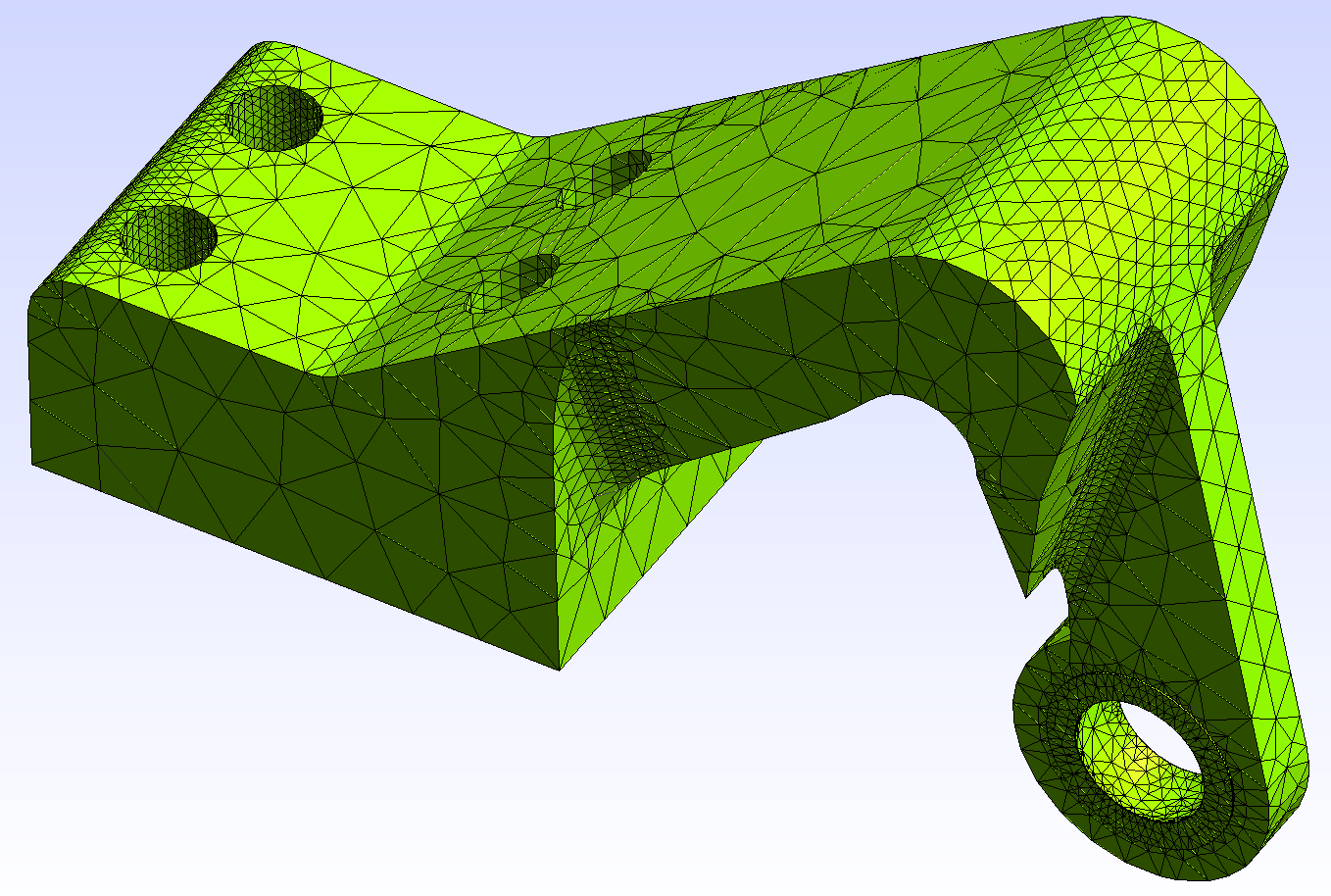}
	\caption{Unseen discretized geometry for PTFONetV3 training. $P^{1}$ conforming tetrahedra finite elements are used for discretization. The mesh consists of $9,223$ nodes and $33,733$ elements} \label{fig:new_geometry}
\end{figure}

We now compare the number of iterations needed for convergence of different preconditioners when a CG iterative solver is used. Figure \ref{fig:num_iterations_several_geometries} shows the number of iterations to convergence of the same combination of preconditioners as in Fig. \ref{fig:num_iterations_single_geometry} plus the addition of NSPOD from a PTFONet trained on a completely new geometry, and hence this PTFONet (PTFONetV3 from now onwards) has not seen any of the geometries in Fig. \ref{fig:geometries} (this new geometry is shown in Fig. \ref{fig:new_geometry}); since PTFONetV3 has only been trained on one geometry, it retains the parameters described in Table \ref{tab:channels} and the weight decay used during its training is $5 \times 10^{-5}$. Figure \ref{fig:num_iterations_several_geometries} shows that the performance of the NSPOD preconditioner using PTFONetV2 also outperforms the \textit{state-of-the-art} preconditioners mentioned in Section \ref{sec:single_geometry}. Even the performance of NSPOD when using PTFONetV3, which was trained on a completely unseen geometry, is comparable to that of BoomerAMG or PETSc's SA-AMG when $m = 1,024$. For Geometry 1, we have that for $m = 1,024$ the number of iterations for NSPOD with PTFONetV2 and PTFONetV3 are 169 and 216, respectively, whilst for BoomerAMG is 218 and for PETSc's SA-AMG is 237. For Geometry 2, we have that for $m = 1,024$ the number of iterations for NSPOD with PTFONetV2 and PTFONetV3 are 85 and 189, respectively, whilst for BoomerAMG is 238 and for PETSc's SA-AMG is 229. For Geometry 3 we have a different trend, where the performance of both versions of NSPOD is poorer than that of the standard AMG preconditioners evaluated; we have that for $m = 1,024$ the number of iterations for NSPOD with PTFONetV2 and PTFONetV3 are 51 and 77, respectively, whilst for BoomerAMG is 33 and for PETSc's SA-AMG is 22.

%%%%%%%%%%%%%%%%%%%%%%%%%%%%%%%%%%%%%%%%%%%

\subsection{Ablation study} \label{sec:ablation_study}
\noindent We also trained a Graph-based version of PTFONet, GraphONet from now onwards, on the geometries shown in Fig. \ref{fig:geometries} at the same time in order to compare its performance against PTFONetV2. Graph-based networks have been used in the last few years to approximate solutions of FE systems because they naturally follow the graph structure of an FE mesh. For instance, Pichi et al. \cite{gcn_fem_1} used a graph autoencoder to encode parametrized PDEs; Despande et al. \cite{gcn_fem_2} used a graph UNet for prediction of PDE solutions. The employed architecture of GraphONet is shown in Fig. \ref{fig:gcn_architecture} and the used parameters are described in Table \ref{tab:channels_several}. We employed the same multigeometry dataset as in the previous section. The relative error on the test set achieved by the trained GraphONet is $32.21\%$.

\begin{figure}
	\centering
	\includegraphics[width=1\linewidth, clip]{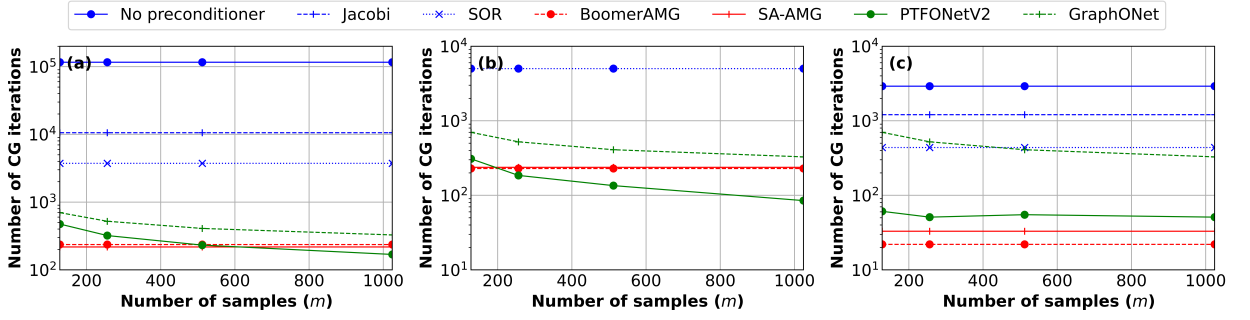}
	\caption{Number of CG iterations \textit{vs} number of samples employed to create the NSPOD preconditioner for the three considered geometries; we employed PTFONetV2 and GraphONet in this scenario. We consider $m = 128, 256, 512, \text{ and } 1024$.
    \textbf{(a)} Geometry 1. 
    \textbf{(b)} Geometry 2, where CG with no preconditioner and/or with Jacobi preconditioner did not converge, thus the corresponding results are not shown.
    \textbf{(c)} Geometry 3.} \label{fig:num_iterations_ablation}
\end{figure}

We now compare the number of iterations needed for convergence of different preconditioners when a CG iterative solver is used. Figure \ref{fig:num_iterations_ablation} shows the number of iterations to convergence of the same combination of preconditioners as in Fig. \ref{fig:num_iterations_several_geometries} but by substituting PTFONetV3 for GraphONet. Figure \ref{fig:num_iterations_ablation} shows that the performance of the NSPOD preconditioner using PTFONetV2 outperforms the NSPOD preconditioner using GraphONet. We have that for $m = 1,024$ the number of iterations for NSPOD with GraphONet for geometries 1, 2, and 3 are, respectively, 328, 114, and 46. Thus, the performance of the NSPOD preconditioner with GraphONet does not improve on the performance of either BoomerAMG or PETSc's SA-AMG for any of the geometries.

%%%%%%%%%%%%%%%%%%%%%%%%%%%%%%%%%%%%%%%%%%%

\subsection{Computational cost} \label{sec:computational_cost}
\noindent The computational cost of the NSPOD preconditioner consists of \textit{offline} and \textit{online} phases. The \textit{offline} cost is dominated by the training of the neural operator (e.g., PTFONet or GraphONet), which is a one-time investment, performed only once prior to the simulation stage.

The \textit{online} cost consists of two components: \textit{assembly} and \textit{application}. During the \textit{assembly} phase, we sample $m \ll n_{\text{nod}}$ latent vectors from a prescribed distribution and pass them through PTFONet. The cost of the PointTransformer pooling layers consists of the cost of MLPs, FPS and $k$NN operations. The cost of MLPs is $\mathcal{O} \left( W \right)$, where $W$ denote the number of weights in a layer. In the case of FPS, it is well known that the cost is roughly $\mathcal{O} \left( r n_{\text{nod}} \log n_{\text{nod}} \right)$, where $r$ represents the pooling ratio. Note that we use a value of $r = 0.25$. The cost of $k$NN is roughly $\mathcal{O} \left( k \left( 1 + r \right) n_{\text{nod}} \log n_{\text{nod}} \right)$, which is a summation of the costs for node feature interpolation and edge adjacency calculation; $k$ is the number of chosen neighbors, which is $k = 16$ in this study. Since the number of weights remains constant regardless of the problem size $n_{\text{nod}}$, the assembly cost is dominated by $\mathcal{O} \left( n_{\text{nod}} \log n_{\text{nod}} \right)$. The computational overhead of the assembly process increases with the problem size; however, it maintains a near-linear complexity rather than exhibiting exponential growth. 
On the other hand, the implementation of GPU-accelerated pooling \cite{johnson2019billion} effectively mitigates the theoretical $\mathcal{O} \left( n_{\text{nod}} \log n_{\text{nod}} \right)$ bottleneck, making it highly scalable for large-scale datasets, achieving near-constant empirical execution time for typical dataset sizes and ensuring high scalability.

Furthermore, the use of GPU parallel computing allows for asynchronous execution with respect to the $m$ samples, providing significant potential for further acceleration. Our numerical experiments indicate that while increasing $m$ enhances the preconditioner's performance, it also adds to the cost; thus, a "sweet spot" for $m$ exists where both acceleration and computational efficiency are optimized.

The cost of applying the preconditioner per iteration, the \textit{application} phase, is comparable to that of an algebraic two-grid method. Let $n_{r}$ be the smoothing step number.
Note that the overall cost for the smoothing step is roughly $\mathcal{O} \left( 2 n_{r} \left( n_{\text{nod}} n_\text{dim} \right) \right)$ \cite{Briggs} since FE matrices are sparse. The combined cost of restriction, prolongation, and the coarse-level solver is $\mathcal{O} \left( m \left( n_{\text{nod}} n_\text{dim} \right) + m^{2} \right)$. Consequently, the total cost per iteration, including smoothing, is $\mathcal{O} \left( 2 n_{r} \left( n_{\text{nod}} n_\text{dim} \right) + m \left( n_{\text{nod}} n_\text{dim} \right) + m^{2} \right)$.

%%%%%%%%%%%%%%%%%%%%%%%%%%%%%%%%%%%%%%%%%%%

\section{Discussion and Conclusion} \label{sec:conclusion}
\noindent In this paper we proposed PTFONet and NSPOD. PTFONet is a PointTransformer-based DeepONet architecture which can effectively predict displacements for linearized solid mechanics on complex CAD geometries. NSPOD is a hybrid preconditioner for Krylov-based linear iterative solvers which uses PTFONet to generate the prolongation operator for a multigrid $V \left( 1, 1 \right)$ structure. We designed a number of numerical experiments involving real CAD geometries, which allow for a comparison of NSPOD's performance against \textit{state-of-the-art} preconditioners for linearized solid mechanics. We also performed an ablation study to assess whether the PointTransformer-based architecture of PTFONet is indeed superior to a graph-based alternative. Finally, we described the computational cost of NSPOD for both its \textit{offline} and \textit{online} stages.

The numerical experiments in Section \ref{sec:single_geometry} show that PTFONet achieves high accuracy when trained on a single geometry for different values of Young's modulus, body forces, and locations of homogeneous Dirichlet BCs. Furthermore, the number of CG iterations achieved by NSPOD with PTFONetV1 is significantly smaller than those from standard preconditioners, Jacobi and SOR; but also smaller than \textit{state-of-the-art} multigrid preconditioners, BoomerAMG and PETSc's SA-AMG, except for Geometry 3. This might be because PTFONet trained on Geometry 3 had the lowest accuracy, but also because Geometry 3 is dominated by ``thin'' features; it is possible that PTFONet only managed to capture a very narrow slice of the solution spectrum due to spectral bias \cite{Rahaman2019}, which might not be enough to warrant a substantial reduction in CG iterations when using the associated NSPOD preconditioner.

The numerical experiments in Section \ref{sec:several_geometries} show that PTFONet also achieves high accuracy when trained on several geometries. This is also performed for different values of Young's modulus, body forces, and locations of homogeneous Dirichlet BCs, within each geometry. Moreover, the number of CG iterations achieved by NSPOD with PTFONetV2 is significantly smaller than those from standard preconditioners, Jacobi and SOR; but also smaller than \textit{state-of-the-art} multigrid preconditioners, BoomerAMG and PETSc's SA-AMG, again with the exception of Geometry 3. The reasoning here is similar to that of PTFONet trained on a single geometry, with spectral bias possibly playing a role on not effectively capturing the portion of the solution spectrum which might be relevant in constructing the prolongation operator. NSPOD with PTFONetV3 achieves a comparable number of CG iterations compared to BoomerAMG and PETSc's SA-AMG, but not significantly lower. This is expected, as PTFONetV3 was trained on a completely different geometry; however, seeing how it still performs comparably to other \textit{state-of-the-art} multigrid preconditioners suggests that there is a degree of geometric transferability.

Our ablation study shows how the performance of NSPOD with a graph-based DeepONet, namely GraphONet, achieves a much lower reduction in CG iterations compared to NSPOD with PTFONetV2. This suggests that a PointTransformer-based DeepONet is indeed superior to an homologous graph-based architecture. However, it is relevant to mention that the accuracy of GraphONet is substantially lower to that of PTFONetV2, i.e. $32.21\%$ vs $16.91\%$ error, which might be the reason why such a relatively low reduction in CG iterations is achieved. Increasing the capacity of GraphONet and retraining with a larger dataset might be a partial solution to this low accuracy. Since FE meshes might have long distance connections between nodes, it is possible that graph-based architectures suffer from oversmoothing \cite{gnn_oversmoothing_1, gnn_oversmoothing_2}, which might lead to a drop in performance. Likewise, whilst connectivity in an FE mesh is extremely important, it is true that connectivity alone does not explain spatial variability of input functions and/or geometric features; spatial coordinates are required to do so. Hence they need to be somehow supplied to all the branches in a DeepONet architecture, which is what occurs in PTFONet. 

This study comes with a number of limitations. We keep variability as constrained as possible; we used spatially constant Young's modulus and body forces, and we only used homogeneous Dirichlet BCs. We only used a few CAD geometries because of the complexity in creating the datasets for each of them; more geometries are required in order to assess the geometric transferability of PTFONet more thoroughly. Also related to this, we only used geometries with $n_\text{nod} \approx 10,000$; much smaller/larger geometries might need modification of the architecture, either changing the pooling ratio or changing the number of pooling operations; these are respectively $0.25$ and $3$ in this study (Fig. \ref{fig:ptfonet_architecture}). The presence of $k$NN and FPS operations in the pooling layers means that large meshes might imply a substantial increase in computational cost even during inference; however, this is considerably reduced by using efficient implementations with KDTrees optimized for GPU. The \textit{offline} stage of NSPOD, i.e. training, is computationally expensive; for instance, PTFONetV2, which has a large capacity, requires around 800 hours of training in a single NVIDIA H100 GPU. Last but not least, the results show that Geometry 3 does not benefit from a large reduction in CG iterations when using NSPOD; complex geometries and/or with ``thin'' features might not benefit from substantial reductions in CG iterations compared to other multigrid preconditioners, which could be due to the presence of spectral bias in PTFONet.

We are also considering a number of ways in which this study could be improved in the future. A first logical step is to extend the study to consider higher variability in terms of spatially-varying input functions, non-homogeneous Dirichlet/Neumann BCs, and number of considered geometries. We should also consider a number of meshes with different numbers of nodes so that we can tune the pooling methodology in the branch network to achieve higher generality. Geometry 3, which comprises of ``thin'' features mainly, benefits the least from a reduction in CG iterations when using NSPOD; this effect might be partially alleviated by increasing the capacity of the architecture, increasing the number of training samples, and/or implementing qualitative changes of the architecture, e.g. such as by including physics-related losses in the loss function. NSPOD is not restricted to a specific architecture but in terms of its practical use in general purpose FEM software it would be advantageous to use architectures which possess a high degree of geometric transferability; whilst geometric transferability is still an open issue in physics ML, PTFONet seems to be allow for a certain degree of geometric transferability when using NSPOD. The latter seems to suggest that hybridizing ML methods with traditional numerical solvers for PDEs is a very powerful tool that might, at least partially, alleviate some of the drawbacks of pure data-based methods.

%%%%%%%%%%%%%%%%%%%%%%%%%%%%%%%%%%%%%%%%%%%

%\appendix

%%%%%%%%%%%%%%%%%%%%%%%%%%%%%%%%%%%%%%%%%%%

% \section{My Appendix}

%%%%%%%%%%%%%%%%%%%%%%%%%%%%%%%%%%%%%%%%%%%

\printcredits

\bibliographystyle{cas-model2-names}

\bibliography{cas-refs}

\end{document}